\documentclass[a4paper,11pt]{article}

\usepackage{amsmath, amssymb, amsfonts, amsthm, enumerate, hyperref, todonotes}

\usepackage[a4paper,top=3cm,bottom=3cm,left=3cm,right=3cm,marginparwidth=3cm]{geometry} 
\usepackage{array}
\usepackage{tasks}
\usepackage{extarrows} 
\usepackage{multirow}
\usepackage{mathtools}
\newtheorem{theorem}{\bf Theorem}[section]
\newtheorem{claim}[theorem]{\bf Claim}
\newtheorem{lemma}[theorem]{\bf Lemma}

\newtheorem{conj}[theorem]{\bf Conjecture}
\newtheorem{prop}[theorem]{\bf Proposition}

\newtheorem{corollary}[theorem]{\bf Corollary}

\newtheorem{obs}[theorem]{\bf Observation}

\theoremstyle{definition}
\newtheorem{defn}[theorem]{\bf Definition}

\newcolumntype{P}[1]{>{\centering\arraybackslash}p{#1}}

\newcommand{\secref}[1]{\hyperref[#1]{Section~\ref*{#1}}}
\newcommand{\thmref}[1]{\hyperref[#1]{Theorem~\ref*{#1}}}
\newcommand{\lemref}[1]{\hyperref[#1]{Lemma~\ref*{#1}}}
\newcommand{\propref}[1]{\hyperref[#1]{Proposition~\ref*{#1}}}
\newcommand{\corref}[1]{\hyperref[#1]{Corollary~\ref*{#1}}}
\newcommand{\defref}[1]{\hyperref[#1]{Definition~\ref*{#1}}}
\newcommand{\obsref}[1]{\hyperref[#1]{Observation~\ref*{#1}}}
\newcommand{\conjref}[1]{\hyperref[#1]{Conjecture~\ref*{#1}}}
\newcommand{\claref}[1]{\hyperref[#1]{Claim~\ref*{#1}}}
\newcommand{\figref}[1]{\hyperref[#1]{Figure~\ref*{#1}}}
\let\eps=\varepsilon

\newcommand{\CC}[1]{\overrightarrow{C\hspace{4pt}}_{\hspace{-5pt}#1}}
\newcommand{\WW}[1]{\overrightarrow{W\hspace{4pt}}_{\hspace{-5pt}#1}}

\let\Olfloor\lfloor
\let\Orfloor\rfloor
\let\Olceil\lceil
\let\Orceil\rceil
\renewcommand{\lfloor}{\left\Olfloor}
\renewcommand{\rfloor}{\right\Orfloor}
\renewcommand{\lceil}{\left\Olceil}
\renewcommand{\rceil}{\right\Orceil}

\usepackage{comment}

\usepackage{tikz}
\usetikzlibrary{arrows.meta,backgrounds,shapes.geometric}
\tikzset{
  blob/.style={ellipse, draw=black!55, fill=black!6, minimum width=12mm, minimum height=7mm, inner sep=0pt, font=\scriptsize},
  cyc/.style={-{Stealth[length=6pt,width=7pt]}, line width=2pt, black!28, shorten <=1pt, shorten >=1pt},
  vtx/.style={circle, fill=black, inner sep=0pt, minimum size=4pt},
  hvtx/.style={circle, draw=black, densely dashed, fill=white, inner sep=0pt, minimum size=5pt},
  arc/.style={-{Stealth[length=4.5pt,width=3.5pt]}, line width=0.55pt, shorten >=0.5pt},
  harc/.style={arc, densely dashed},
  vlab/.style={font=\footnotesize, inner sep=1.5pt},
  pcap/.style={font=\small},
  hl/.style={line width=6.5pt, line cap=round, line join=round, orange!45},
  hlb/.style={hl, cyan!40},
}
\newcommand{\blobsfour}{%
  \path[use as bounding box] (-2.35,-2.1) rectangle (2.35,2.1);
  \node[blob] (B0) at (90:1.7) {$V_{p(v)}$};
  \node[blob] (B1) at (0:1.7) {$V_{p(v)+1}$};
  \node[blob] (Bm2) at (270:1.7) {$V_{p(v)-2}$};
  \node[blob] (Bm1) at (180:1.7) {$V_{p(v)-1}$};
  \foreach \a/\b in {Bm1/B0, B0/B1, B1/Bm2, Bm2/Bm1} \draw[cyc] (\a) -- (\b);}
\newcommand{\blobsfive}{%
  \path[use as bounding box] (-2.41,-2.6) rectangle (2.41,2.25);
  \node[blob] (B0) at (90:1.85) {$V_{p(v)}$};
  \node[blob] (B1) at (18:1.85) {$V_{p(v)+1}$};
  \node[blob] (B2) at (-54:1.85) {$V_{p(v)+2}$};
  \node[blob] (Bm2) at (-126:1.85) {$V_{p(v)-2}$};
  \node[blob] (Bm1) at (162:1.85) {$V_{p(v)-1}$};
  \foreach \a/\b in {Bm1/B0, B0/B1, B1/B2, B2/Bm2, Bm2/Bm1} \draw[cyc] (\a) -- (\b);}

\begin{document}
\title{
Minimum outdegree forcing directed cycles\thanks{This work was supported by the National Science Centre grant 2021/42/E/ST1/00193.}
}

\date{}

\author{
Andrzej Grzesik\thanks{Faculty of Mathematics and Computer Science, Jagiellonian University, {\L}ojasiewicza 6, 30-348 Krak\'{o}w, Poland. E-mail: {\tt andrzej.grzesik@uj.edu.pl}.}
\and 
Justyna Jaworska\thanks{Faculty of Mathematics and Computer Science, Jagiellonian University, {\L}ojasiewicza 6, 30-348 Krak\'{o}w, Poland. E-mail: {\tt justynajoanna.jaworska@student.uj.edu.pl}.}
\and 
Piotr Kuc\thanks{Faculty of Mathematics and Computer Science, Jagiellonian University, {\L}ojasiewicza 6, 30-348 Krak\'{o}w, Poland. E-mail: {\tt piotr.kuc@student.uj.edu.pl}.}
}

\maketitle

\begin{abstract}
The Caccetta--Häggkvist conjecture, one of the best known open problems on oriented graphs, asserts that every oriented graph on $n$ vertices with minimum outdegree at least $n/\ell$ contains a directed cycle of length at most $\ell$. Motivated by it, Kelly, Kühn and Osthus asked for the minimum degree condition that forces a directed cycle of length exactly $\ell$. The semidegree version of this problem was recently solved, but the methods used there rely on bounding the directed diameter, which requires control over the indegrees. In the setting of the Caccetta--Häggkvist conjecture, where only the outdegrees are bounded, the answer was known only for $\ell$ not divisible by~$3$.

We determine the optimal minimum outdegree forcing a directed cycle of length $\ell$ for every sufficiently large $\ell$ and $n$. This establishes the outdegree version of the Kelly--Kühn--Osthus conjecture for all large $\ell$, up to a small correction of the additive term, which in general cannot be avoided.
%
\end{abstract}

\section{Introduction}

One of the best known and most studied open problems concerning oriented graphs is the following conjecture from 1978.

\begin{conj}[Caccetta, Häggkvist \cite{caccetta_haggkvist}]\label{con:CH}
Every oriented graph $G$ on $n$ vertices with minimum outdegree $\delta^+(G) \geq \frac{n}{\ell}$ contains a directed cycle of length at most $\ell$.
\end{conj}

Here, by an \emph{oriented graph} we understand a directed graph without loops and symmetric pairs of directed arcs, that is, an orientation of a simple graph. Furthermore,
$\delta^+(G)$ and $\delta^-(G)$ denote the \emph{minimum outdegree} and the \emph{minimum indegree} of an oriented graph $G$, respectively.

Numerous partial results towards Conjecture~\ref{con:CH} are known. These include proofs with an additive error term in the bound on the cycle length, proofs in which the minimum outdegree assumption is strengthened by a multiplicative constant, and results obtained under additional forbidden-subgraph assumptions. For these and related developments concerning the Caccetta--Häggkvist conjecture, we refer the reader to the survey~\cite{CH_survey}.

Even if we strengthen the assumption of Conjecture~\ref{con:CH} by imposing a lower bound on the \emph{minimum semidegree} $\delta^{\pm}(G) = \min(\delta^+(G), \delta^-(G))$ instead of the minimum outdegree, the essential difficulty seems to persist and the resulting weaker statement is open as well, except the case $\ell=3$ recently claimed to be solved by OpenAI \cite{OpenAI}.

Motivated by this conjecture, Kelly, Kühn and Osthus~\cite{kko} proposed in 2010 the following problem: for a given integer $\ell \ge 4$, determine the minimum semidegree that forces an oriented graph to contain a directed cycle of length \emph{exactly}~$\ell$. Despite the similar formulation, this problem is substantially different, and the answer depends on the divisibility properties of $\ell$. Indeed, Kelly, Kühn and Osthus conjectured the following.

\begin{conj}[Kelly, Kühn, Osthus \cite{kko}]\label{con:KKO}
Let $\ell \geq 4$ be an integer and let $k$ be the smallest integer greater than $2$ that does not divide $\ell$. Then there exists $n_0$ such that every oriented graph $G$ on $n \geq n_0$ vertices with minimum semidegree $\delta^{\pm}(G) \geq \lfloor \frac nk \rfloor + 1$ contains a directed cycle of length $\ell$.
\end{conj}

Kelly, Kühn and Osthus~\cite{kko} proved Conjecture~\ref{con:KKO} in the case $k=3$, that is, for cycle lengths not divisible by $3$. They also showed an asymptotic statement (with the semidegree bound $(1+\eps)\frac{n}{k}$) in the cases $k=4$ and $k=5$. Later, Kühn, Osthus and Piguet~\cite{kop} proved a similar asymptotic statement for the remaining values of $k$, under the additional assumption that $\ell \geq 10^7 k^6$.
Conjecture~\ref{con:KKO} was then resolved in 2023, up to a slight adjustment of the semidegree bound.

\begin{theorem}[Grzesik, Volec \cite{grzesik2023degree}]\label{thm:delta_pm}
    Let $\ell \geq 4$ be an integer and let $k$ be the smallest integer greater than $2$ that does not divide $\ell$. Then there exists $n_0$ such that for every oriented graph $G$ on $n \geq n_0$ vertices:
    \begin{itemize}
        \item if $\ell \not\equiv 3 \bmod 12$ and $\delta^\pm (G) \geq \frac nk + \frac{k-1}{2k}$, then $G$ contains a directed cycle of length $\ell$,
        \item if $\ell \equiv 3 \bmod 12$ and $\delta^\pm (G) \geq \frac n4 + \frac 14$, then $G$ contains a directed cycle of length $\ell$.
    \end{itemize}
    Moreover, the above bounds are best possible.
\end{theorem}

As with the Caccetta--Häggkvist conjecture, Conjecture~\ref{con:KKO} seems to hold in a stronger form, where the minimum semidegree is replaced by the minimum outdegree. However, the techniques used in the aforementioned results cannot be easily adapted to this stronger setting, as they heavily rely on bounding the directed diameter, which cannot be done without an assumption on the indegrees.

As a first step towards the outdegree version of Conjecture~\ref{con:KKO}, Czygrinow, Molla, Nagle and Oursler~\cite{even_rainbow} proved that for every integer $\ell \geq 4$, there exists $n_0$ such that every oriented graph $G$ on $n \geq n_0$ vertices with minimum outdegree $\delta^+(G) \geq \frac{n}{3} + \frac{1}{3}$ contains a directed cycle of length $\ell$. In particular, the outdegree version of Conjecture~\ref{con:KKO} holds whenever $\ell$ is not divisible by $3$.

In this paper, we improve the above result for large $\ell$ divisible by $3$, providing an optimal outdegree bound. In particular, this shows that the outdegree version of Conjecture~\ref{con:KKO} holds for every sufficiently large $\ell$, up to the same kind of slight adjustment of the degree bound as in Theorem~\ref{thm:delta_pm}.

\begin{theorem}\label{thm:main_thm}
    Let $\ell$ be an integer divisible by $3$ and let $k$ be the smallest integer greater than $2$ that does not divide $\ell$. If $\ell > 20k^2$, then there exists $n_0$ such that every oriented graph $G$ on $n \geq n_0$ vertices with $\delta^+(G) \geq \frac nk + \frac{k-1}{k}$ contains a directed cycle of length $\ell$.
    Moreover, this outdegree bound is best possible if $\ell \not\equiv 3,6 \bmod 12$ and $\ell \not\equiv 24,48 \bmod 60$.
\end{theorem}

Note that $k$ is always a prime power, and so $k > 9$ forces $k \geq 11$. In particular, when $k > 9$, the assumption $\ell > 20k^2$ of Theorem~\ref{thm:main_thm} is redundant, as every such $\ell$ satisfies it: all of the integers $3, 4, \dots, k-1$ divide $\ell$, and hence $\ell \geq \operatorname{lcm}(3, 4, \dots, k-1) > 20k^2$. Thus, the assumption $\ell > 20k^2$ can be replaced by $\ell > 1620$. A more careful analysis of our proof would reduce this bound even further, to $\ell > 180$, if needed.

Observe that the outdegree bound in Theorem~\ref{thm:main_thm} is slightly larger than the semidegree bound in Theorem~\ref{thm:delta_pm}. At an intuitive level, this is because a larger family of extremal constructions is available when only the outdegree condition has to be maintained.

However, similarly as with $\ell \equiv 3 \bmod 12$ in Theorem~\ref{thm:delta_pm}, there are cases in which the best known construction does not match the general bound in Theorem~\ref{thm:main_thm}. It turns out that in these cases the outdegree bound can indeed be improved, and the following theorem determines its optimal value.

\begin{theorem}\label{thm:special_cases}
    Let $\ell$ be an integer divisible by $3$ and let $k$ be the smallest integer greater than $2$ that does not divide $\ell$. Assume that $\ell > 20k^2$, $G$ is an oriented graph on $n$ vertices, with $n$ large enough, and one of the following assumptions is satisfied:
    \begin{itemize}
        \item $\ell \equiv 3 \bmod 12$ and $\delta^+ (G) \geq \frac n4 + \frac 14$;
        \item $\ell \equiv 6 \bmod 12$ and $\delta^+ (G) \geq \frac n4 + \frac 24$;
        \item $\ell \equiv 24 \bmod 60$ and $\delta^+ (G) \geq \frac n5 + \frac 35$;
        \item $\ell \equiv 48 \bmod 60$ and $\delta^+ (G) \geq \frac n5 + \frac 25$.
    \end{itemize}
    Then $G$ contains a directed cycle of length $\ell$. Moreover, these bounds are best possible.
\end{theorem}

The rest of the paper is organized as follows. In Section~\ref{section:prelims}, we introduce our notation and several lemmas we will use in the subsequent proofs. In Section~\ref{section:constructions}, we present the extremal constructions for all values of $\ell$. Sections~\ref{section:main_thm_proof} and \ref{section:special_cases} contain the proofs of Theorems~\ref{thm:main_thm} and \ref{thm:special_cases}, respectively.

\section{Notation and Preliminaries}\label{section:prelims}

\subsection{Notation}

A \emph{directed graph} $G$ is a pair $(V(G), E(G))$, where $E(G) \subseteq V(G) \times V(G)$. Elements of $E(G)$ are called \emph{arcs}. When $(u, v) \in E(G)$, we say there is an arc \emph{from} $u$ \emph{to} $v$. We often write $uv$ instead of $(u, v)$ for brevity. In all directed graphs we consider, \emph{loops}, i.e., arcs of the form $vv$, are not allowed.

\emph{Oriented graphs} are directed graphs with an additional restriction -- for $u, v \in V(G)$, if $uv \in E(G)$, then $vu \not\in E(G)$ (this in particular implies there are no loops). As we are mostly concerned with oriented graphs, we refer to them simply as \emph{graphs} whenever this is clear from the context.

In a directed graph, an \emph{outneighbor} of a vertex $v$ is a vertex $u$ such that $vu$ is an arc. Similarly, $u$ is an \emph{inneighbor} of $v$ if $uv$ is an arc.

The \emph{outneighborhood} $N^+(v)$ of a vertex $v$ is the set of its outneighbors. Symmetrically, the \emph{inneighborhood} $N^-(v)$ is the set of the inneighbors of $v$. The \emph{outdegree} and \emph{indegree} of $v$ are $d^+(v) := |N^+(v)|$ and $d^-(v) := |N^-(v)|$, respectively. Whenever the underlying graph is not clear from the context, we indicate it by a subscript, writing $N^+_G(v)$, $d^+_G(v)$, and so on. Define:
\begin{align*}
&\delta^+(G) = \min_{v \in V(G)} d^+(v)
&&\Delta^+(G) = \max_{v \in V(G)} d^+(v) &\ &\text{minimum/maximum outdegree} \\
&\delta^-(G) = \min_{v \in V(G)} d^-(v)
&&\Delta^-(G) = \max_{v \in V(G)} d^-(v) &\  &\text{minimum/maximum indegree} \\
&\delta^{\pm}(G) = \min(\delta^+(G), \delta^-(G)) &&&\ &\text{minimum semidegree}
\end{align*}

We extend the notion of out- and inneighborhood to sets of vertices -- for $A \subseteq V(G)$, we set $N^+(A) := \bigcup_{a \in A} N^+(a)$, and define $N^-(A)$ analogously.

For $A,B \subseteq V(G)$, we set $E(A, B) := E(G) \cap (A \times B)$ (the set of arcs from $A$ to $B$ in the graph $G$) and $e(A, B) := |E(A, B)|$.

A \emph{blow-up} of a graph $G$ on $k$ vertices $v_1, v_2, \ldots, v_k$ is any graph $G'$ whose vertex set is partitioned into sets $V_1, V_2, \ldots, V_k$ and in which, for all $i, j$ and all $u \in V_i$, $v \in V_j$, there is an arc $uv$ if and only if there is an arc $v_iv_j$ in $G$. The sets $V_i$ are called \emph{blobs}; note that they are independent sets, as $G$ has no loops. We sometimes refer to a blow-up of a graph $G$ as a \emph{$G$-blow-up}.

A blow-up is \emph{balanced} if $\bigl| |V_i| - |V_j| \bigr| \leq 1$ for any $i, j$. Note that for every integer $n \geq k$ there exists a balanced blow-up of $G$ on $n$ vertices. By $G[t]$ we denote the blow-up of $G$ in which $|V_1| = |V_2| = \ldots = |V_k| = t$.

A graph $H = (V_H, E_H)$ is a \emph{subgraph} of $G = (V_G, E_G)$ if $V_H \subseteq V_G$ and $E_H \subseteq E_G$. The subgraph is \emph{spanning} if $V_H = V_G$. Moreover, we say that $H$ is an \emph{induced subgraph} of $G$ if $E_H = E_G \cap (V_H \times V_H)$. The induced subgraph of $G$ on a vertex set $V_H$ is denoted by $G[V_H]$.

For a graph $F$, we say $G$ is \emph{$F$-free} if $G$ does not contain $F$ as a subgraph. For a family $\mathcal{F}$ of graphs, we say that $G$ is \emph{$\mathcal{F}$-free} if $G$ does not contain any element of $\mathcal{F}$ as a subgraph.

A graph $H$ is \emph{homomorphic} to a graph $F$ if there exists a function $\phi : V(H) \rightarrow V(F)$ such that $\phi(h_1)\phi(h_2) \in E(F)$ for every arc $h_1h_2 \in E(H)$; such a function is called a \emph{homomorphism} from $H$ to $F$. Moreover, $F$ is a \emph{homomorphic image} of $H$ if there exists a homomorphism $\phi$ from $H$ to $F$ which is surjective and under which every arc of $F$ is the image of an arc of $H$; formally, if for every $f_1f_2 \in E(F)$ there is an arc $h_1h_2 \in E(H)$ with $f_1 = \phi(h_1)$ and $f_2 = \phi(h_2)$.

Observe that $H$ is a subgraph of $F[t]$ for some $t$ if and only if $F$ contains a homomorphic image of $H$.

For an integer $\ell \geq 3$, by $\CC{\ell}$ we denote the \emph{directed cycle} of length $\ell$, that is, the oriented graph with vertices $v_1, v_2, \ldots, v_{\ell}$ and arcs $v_1v_2, v_2v_3, \ldots, v_{\ell}v_1$. We also refer to $\CC{\ell}$ as the \emph{$\ell$-cycle}.

A \emph{directed walk} of length $\ell$ is a sequence of not necessarily distinct vertices $v_0, v_1, \ldots, v_{\ell}$ such that there are arcs $v_0v_1, v_1v_2, \ldots, v_{\ell-1} v_{\ell}$. A directed walk is \emph{closed} if $v_0 = v_{\ell}$. By $\WW{\ell}$ we denote the family of homomorphic images of $\CC{\ell}$. Note that a graph $G$ is $\WW{\ell}$-free if and only if it does not contain any closed walk of length $\ell$.

A \emph{directed path} of length $\ell$ is an oriented graph $P$ on $\ell+1$ vertices $v_0, v_1 \ldots, v_{\ell}$ with $\ell$ arcs $v_0v_1, v_1v_2, \ldots, v_{\ell-1}v_{\ell}$. We say that such a path \emph{starts} in $v_0$ and \emph{ends} in $v_\ell$. The vertices $v_0, v_{\ell}$ are called the \emph{end vertices} of $P$, whereas $v_1, v_2, \ldots, v_{\ell-1}$ are called the \emph{internal vertices}.

The \emph{distance} from a vertex $a$ to a vertex $b$ in an oriented graph $G$ is the length of a shortest directed path in $G$ which starts in $a$ and ends in $b$; if no such path exists, the distance is infinite.

\subsection{Frobenius Coin Problem}

The Frobenius Coin Problem asks, given a sequence of positive integers $(a_1, a_2, \ldots, a_m)$ (``coin denominations''), for the largest integer $N$ which is not expressible in the form $x_1a_1 + x_2a_2 + \ldots + x_ma_m$ for some nonnegative integers $x_1, x_2, \ldots, x_m$. Note that, in order for the problem to be well defined, $\gcd(a_1, a_2, \ldots, a_m) = 1$ must hold.

We use the following classic result by Sylvester, which solves the $m=2$ case.

\begin{prop}[\cite{sylvester}]\label{prop:sylvester}
Let $a, b$ be coprime positive integers. If a positive integer $\ell$ satisfies ${\ell \geq (a-1)(b-1)}$, then $\ell$ can be represented as $xa + yb$, where $x, y$ are nonnegative integers.
\end{prop}

In our applications we need the following, slightly generalized, form.

\begin{corollary}\label{cor:frobenius}
Let $a, b, \ell$ be positive integers such that $\gcd(a, b) \mid \ell$. If $\ell \geq \frac{(a-1)(b-1)}{\gcd(a,b)}$, then $\ell$ can be represented as $xa + yb$, where $x, y$ are nonnegative integers.
\end{corollary}

\begin{proof}
Let $d = \gcd(a, b)$ and define $a' = \frac ad$, $b' = \frac bd$ and $\ell' = \frac \ell d$. Then $a', b', \ell'$ are positive integers, $\gcd(a', b') = 1$, and
$$d(a'-1)(b'-1) \leq \frac{1}{d}(da'-1)(db'-1) = \frac{1}{d}(a-1)(b-1) \leq \ell,$$
so that $\ell' \geq (a'-1)(b'-1)$. By Proposition~\ref{prop:sylvester}, $\ell' = xa' + yb'$ for some nonnegative integers $x, y$, and thus $\ell = xa + yb$.
\end{proof}

\subsection{Sumsets}

We are going to operate on sets of remainders modulo $k$. For $A, B \subseteq \mathbb{Z}_k$, define
$$A + B := \{ a + b \bmod k : a \in A, b \in B \}.$$

Given a family of sets of remainders $A_1, A_2, \ldots, A_m \subseteq \mathbb{Z}_k$, our aim is to \emph{compose} the special remainder $\ell \bmod k$ using at most one remainder from each of the sets $A_i$. To put this formally, we ask whether $(\ell \bmod k) \in (A_1 \cup \{0\})+\ldots + (A_m \cup \{0\})$. Define
$$\mathcal{C}(A_1, A_2, \ldots, A_m) := ((A_1 \cup \{0\})+\ldots + (A_m \cup \{0\})) \setminus \{ 0 \}$$
to be the set of ``composable'' non-zero remainders modulo $k$.

We quote the following two results from \cite{grzesik2023degree}. Even though they were initially stated for $\ell$ divisible by $12$, generalizing to $\ell$ divisible by $6$ is immediate.

\begin{lemma}[\cite{grzesik2023degree}] \label{lemma:remainders_sets}
Let $\ell$ be a multiple of $6$, let $k$ be the smallest integer that does not divide $\ell$, and let $A_1, \ldots, A_m$ be a collection of $m$ sets of non-zero remainders modulo $k$. If $\ell \bmod k$ does not belong to $\mathcal{C} (A_1, \ldots, A_m)$, then ${|\mathcal{C}(A_1, \ldots, A_m)| \geq |A_1| + \ldots + |A_m|}$.

In particular, if $|A_1| + \ldots + |A_m| \geq k-1$, then $\ell \bmod k \in \mathcal{C}(A_1, \ldots, A_m)$.
\end{lemma}

Considering the special case in which every set $A_i$ has exactly one element leads to the following corollary.
\begin{corollary}[\cite{grzesik2023degree}]\label{cor:remainders_plain}
    Let $\ell$ be a multiple of $6$ and let $k$ be the smallest integer that does not divide $\ell$. If $S$ is a multiset of $k-1$ non-zero remainders modulo $k$, then there exists a sub-multiset $S' \subseteq S$ such that $$ \ell \equiv \sum_{x\in S'} x \bmod{k}.$$

\end{corollary}

\subsection{Szemerédi Regularity Lemma}

Lastly, we will use the following proposition, which is a standard application of the Szemerédi Regularity Lemma for directed graphs~\cite{directed_regularity}.

\begin{prop}\label{prop:graph_cleanup}
Let $\ell \geq 4$ be an integer and let $\eps \in (0, 1)$. There exists $n_0 = n_0(\ell, \eps)$ with the following property. If $G$ is a $\CC{\ell}$-free oriented graph on $n \geq n_0$ vertices, then $G$ contains a $\WW{\ell}$-free subgraph $H$ on at least $(1 - \eps)n$ vertices such that $d^+_H(v) \geq d^+_G(v) - \eps n$ for each $v \in V(H)$.

In particular, if an $n$-vertex oriented graph $G$ is $\CC{\ell}$-free, then one can remove $o(n)$ vertices and $o(n^2)$ arcs from $G$ to obtain a $\WW{\ell}$-free subgraph $H$ with ${\delta^+(H) \geq \delta^+(G) - o(n)}$.
\end{prop}

\section{Constructions}\label{section:constructions}

We now present the best known extremal constructions.
Recall that $k$ denotes the smallest integer greater than $2$ that does not divide $\ell$. Our first aim is to construct arbitrarily large $\CC{\ell}$-free graphs $G$ on $n$ vertices with $\delta^+(G) = \frac{n}{k} + \frac{k-2}{k}$, which shows that the outdegree bound in Theorem~\ref{thm:main_thm} cannot be lowered.

To this end, take $n \equiv 2 \bmod k$ and consider the balanced blow-up $G_0$ of $\CC{k}$ on $n+k-2$ vertices; all its blobs are then of the same size $\frac{n+k-2}{k}$. Since every cycle in a $\CC{k}$-blow-up has length divisible by $k$, and $k$ does not divide $\ell$, the graph $G_0$ contains no $\CC{\ell}$. Moreover, $\delta^+(G_0) = \frac{n+k-2}{k} = \frac nk + \frac{k-2}{k}$.

The goal is to decrease the number of vertices to $n$ while preserving the minimum outdegree and the $\CC{\ell}$-freeness of our construction. We do this iteratively, by building a sequence of $k-1$ graphs, $G_0, G_1, \ldots, G_{k-2}$. Denote
$$R_i = \{ r \in \mathbb{Z}_k : G_i \text{ contains } \CC{m} \text{ for some } m \equiv r \bmod{k} \}.$$
Clearly, $R_0 = \{ 0 \}$. In each step of our construction we ensure that:


\begin{enumerate}
    \item\label{item:degree} $\delta^+(G_{i+1}) = \delta^+(G_{i})$;
    \item\label{item:vertices} $|V(G_{i+1})| = |V(G_{i})| - 1$;
    \item\label{item:remainders} $|R_{i+1}| = |R_i| + 1$.
\end{enumerate}
At the end we will obtain a graph $G_{k-2}$ on $n$ vertices with $\delta^+(G_{k-2}) = \delta^+(G_0) = \frac nk + \frac{k-2}{k}$. We aim to have $R_{k-2} = \mathbb{Z}_k \setminus \{ \ell \bmod{k} \}$, which will make $G_{k-2}$ the desired extremal construction.

\begin{figure}[ht]
    \centering
    \includegraphics[width=0.3\linewidth]{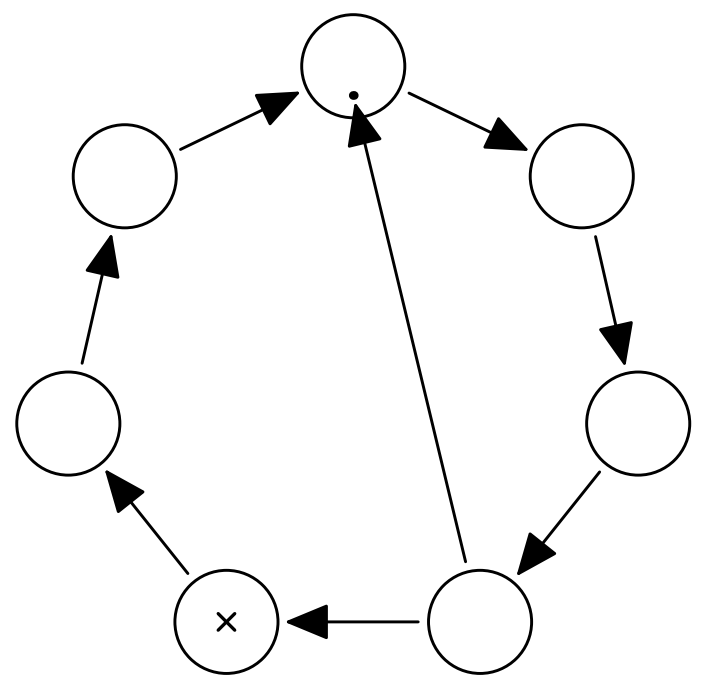}
    \caption{Example of a $4$-maneuver for $k=7$.}
    \label{fig:maneuver}
\end{figure}

Since we only remove vertices and add arcs in our process, each vertex in the graphs $G_i$ has a unique blob of the initial $\CC{k}$-blow-up $G_0$ associated with it. Refer to these blobs as $V_1, \ldots, V_k$, with indices taken modulo $k$. For each $r \in \{3, 4, \ldots, k-1\}$, we define an \emph{$r$-maneuver} as follows:
\begin{itemize}
    \item pick a special vertex $v$, which was not chosen before, from a blob $V_j$ such that no vertex of $V_{j+r-1}$ is an outneighbor of $v$;
    \item add arcs from all vertices in $V_{j+r-1}$ to the special vertex $v$;
    \item delete a vertex (which was not special in any of the previous maneuvers) from $V_{j+r}$.
\end{itemize}
See Figure~\ref{fig:maneuver} for an example. The restriction on the blob $V_j$ guarantees that $G_{i+1}$ is again an oriented graph, that is, that the added arcs do not create a symmetric pair. Since each maneuver excludes at most one of the $k$ blobs, an admissible blob always exists. Note also that the values $r \leq 2$ are not allowed, since the added arcs would create a loop (for $r=1$) or a symmetric pair (for $r=2$).

It is straightforward to see that Properties~\ref{item:degree} and~\ref{item:vertices} are preserved. Indeed, deleting a vertex of $V_{j+r}$ decreases the outdegree of every vertex of $V_{j+r-1}$ by one, which is exactly compensated by the arcs added to $v$, while all the remaining outdegrees are unchanged; here we use that the deleted vertex is not special, and hence its only inneighbors are the vertices of $V_{j+r-1}$. Moreover, if $G_{i+1}$ was obtained from $G_i$ through an $r$-maneuver, then $R_{i+1} = R_i + \{0, r\}$, as we can use the special vertex to make one shortcut and return to the starting blob after $r$ steps rather than $k$. Observe that if in each step we apply either an $r$-maneuver or a $(k-r)$-maneuver for a fixed $r$ coprime with $k$, then Property~\ref{item:remainders} holds.

We will now describe constructions (i.e.~choices of maneuvers) for various $\ell$. We start with the case $k \geq 7$.

For odd $k$, take $r = \frac{k-1}{2}$. Since $2r \equiv -1 \bmod k$, the number $r$ is coprime with $k$, and so there exists $a \in \{1, 2, \ldots, k-1\}$ such that ${ar \equiv \ell \bmod{k}}$; note that $a \neq 0$, as $k$ does not divide $\ell$. Now, perform $(a-1)$ times an $r$-maneuver and $(k-a-1)$ times a $(k-r)$-maneuver. Altogether this makes $k-2$ maneuvers, and since a $(k-r)$-maneuver contributes the remainder $-r$, we end up with
$$R_{k-2} = \{ cr \bmod k : -(k-a-1) \leq c \leq a-1 \} = \mathbb{Z}_k \setminus \{ \ell \bmod{k} \},$$
because $a$ is the only remainder modulo $k$ missing from the above range of $c$.

For even $k$, recall that $k$ is always a prime power, so $k = 2^s$ for some integer $s \geq 3$. Therefore, $r=3$ is coprime with $k$, and we may proceed as above: we take $a$ such that $3a \equiv \ell \bmod{k}$ and apply $(a-1)$ times a $3$-maneuver and $(k-a-1)$ times a $(k-3)$-maneuver.

For $k=5$, only $3$- and $4$-maneuvers are available, and there is no $r$ for which both an $r$-maneuver and a $(k-r)$-maneuver could be used. Therefore, we cannot use the same method as above. Table~\ref{tab:constructions_k5} lists the best possible constructions with the available maneuvers. Here, as in the general construction, we assume that $n$ is such that $k$ divides the number of vertices of the initial blow-up; the same applies to Table~\ref{tab:constructions_k4} below.
\begin{table}[ht]
    \centering
    \renewcommand{\arraystretch}{1.5} 
    \begin{tabular}{|P{0.09\linewidth}|p{0.57\linewidth}|p{0.24\linewidth}|} 
        \hline
        \textbf{$\ell \bmod 5$} & \textbf{Construction} & \textbf{Remarks} \\ 
        \hline
        $1$ & apply three times a $4$-maneuver to a balanced \hbox{$\CC{5}$-blow-up} on $n+3$ vertices & Matches Theorem~\ref{thm:main_thm}\\ 
        \hline
        $2$ & apply three times a $3$-maneuver to a balanced \hbox{$\CC{5}$-blow-up} on $n+3$ vertices & Matches Theorem~\ref{thm:main_thm}\\ 
        \hline
        $3$ & apply once a $4$-maneuver to a balanced $\CC{5}$-blow-up on $n+1$ vertices & Matches Theorem~\ref{thm:special_cases}\\ 
        \hline
        $4$ & apply twice a $3$-maneuver to a balanced $\CC{5}$-blow-up on $n+2$ vertices & Matches Theorem~\ref{thm:special_cases}\\ 
        \hline
    \end{tabular}
    \caption{Best known constructions for $k=5$.}
    \label{tab:constructions_k5}
\end{table}

For $k=4$, we are limited to $3$-maneuvers, resulting in the constructions in Table~\ref{tab:constructions_k4}.

\begin{table}[ht]
    \centering
    \renewcommand{\arraystretch}{1.5} 
    \begin{tabular}{|P{0.09\linewidth}|p{0.57\linewidth}|p{0.24\linewidth}|} 
        \hline
        \textbf{$\ell \bmod 4$} & \textbf{Construction} & \textbf{Remarks} \\ 
        \hline
        $1$ & apply twice a $3$-maneuver to a balanced $\CC{4}$-blow-up on $n+2$ vertices & Matches Theorem~\ref{thm:main_thm} \\ 
        \hline
        $2$ & apply once a $3$-maneuver to a balanced $\CC{4}$-blow-up on $n+1$ vertices & Matches Theorem~\ref{thm:special_cases} \\ 
        \hline
        $3$ & a balanced $\CC{4}$-blow-up on $n$ vertices & Matches Theorem~\ref{thm:special_cases} \\ 
        \hline
    \end{tabular}
    \caption{Best known constructions for $k=4$.}
    \label{tab:constructions_k4}
\end{table}

By Theorems~\ref{thm:main_thm} and \ref{thm:special_cases}, these constructions are best possible for large enough cycles. We believe that they are best possible for small cycles as well.

\section{Proof of Theorem \ref{thm:main_thm}} \label{section:main_thm_proof}

The first main step of the proof of Theorem~\ref{thm:main_thm} is to show that any $\CC{\ell}$-free oriented graph $G$ on $n$ vertices with $\delta^+(G) > \frac nk$ contains a subgraph $H$ homomorphic to $\CC{k}$ having semidegree $\delta^\pm(H) \ge \frac{n}{k}-o(n)$. This reasoning is split into two cases based on the parity of $\ell$ in the next two subsections. The second step is to determine more structure of $G$ by using the outdegree bound on $G$ and the semidegree bound on $H$, in particular proving that one can assume that $H$ is an induced subgraph of $G$. This is done in Subsection~\ref{section:framework}, where we introduce useful notation and prove observations that will be also useful in the later proof of Theorem~\ref{thm:special_cases}. The conclusion of the proof of Theorem~\ref{thm:main_thm} is in Subsection~\ref{sec:main_proof_finish} and contains the only place where the exact bound $\delta^+(G) \geq \frac nk + \frac{k-1}{k}$ is used.

In the proofs that follow, we use the following approximate result towards the Caccetta--Häggkvist conjecture, proved in~\cite{grzesik2023degree}.

\begin{lemma}[\cite{grzesik2023degree}]\label{lemma:girth_3m}
    Let $m$ be a positive integer. Every oriented graph $F$ on $N$ vertices with $\delta^+(F) \geq (1 - 0.6535^{1/m})N$ contains a directed cycle of length at most $3m$. 
\end{lemma}

\subsection{Approximate structure in the even case}\label{sec:even_case_approx}

The objective of this section is proving the following lemma.

\begin{lemma}\label{lemma:stability_even_case}
    Let $\ell$ be an integer divisible by $6$ and denote by $k$ the smallest integer greater than $2$ which does not divide $\ell$. If $\ell > 20k^2$, then every $\CC{\ell}$-free oriented graph $G$ on $n$ vertices with $\delta^+(G) > \frac nk$ contains a subgraph $H$ homomorphic to $\CC{k}$ with $\delta^{\pm}(H) \geq \frac nk - o(n)$.
\end{lemma}

\begin{proof}

Apply Proposition \ref{prop:graph_cleanup} to the graph $G$ to obtain a $\WW{\ell}$-free subgraph $H$ with $\delta^+(H) \ge \frac nk - o(n)$. 
We will show that after removing $o(n)$ vertices from $H$, we obtain the desired graph. 

Note that, by the minimality of $k$, every integer smaller than $k$ divides $\ell$. In particular, if $H$ contained a directed cycle of length $p \leq k-1$, then $p \geq 3$ and $p \mid \ell$, so traversing this cycle $\frac \ell p$ times would give a closed walk of length $\ell$ in $H$, contradicting its $\WW{\ell}$-freeness. We use this observation repeatedly.

We start with the following claim, that resembles \cite[Claim~4.1]{grzesik2023degree} and will allow us to build directed cycles. In contrast to \cite[Claim~4.1]{grzesik2023degree}, we do not have an indegree assumption, thus we cannot claim existence of a short directed path between any two vertices, but only of a short directed path from every vertex to one particularly chosen special vertex. 
\begin{claim}\label{clm:short_paths_general}
    There is a vertex $v_0 \in V(H)$ such that for every $v \in V(H)$, there is a directed path of length at most $10 \lfloor \frac{k-1}{3} \rfloor$ from $v$ to $v_0$.
\end{claim}

\begin{proof}
    Let $m := \lfloor \frac{k-1}{3} \rfloor$, $c := 1 - 0.6535^{1/m}$. Take any vertex $v \in V(H)$. We define $X_1 := N^+(v)$ and $X_{i} := X_{i-1} \cup N^+(X_{i-1})$ for every $i \in \{2,\ldots,5m\}$.

    Observe that for each $i \in \{2, \ldots, 5m \}$, the induced graph $H[X_{i-1}]$ cannot contain a directed cycle of length at most $3m \leq k-1$, as otherwise $H$ would contain a closed walk of length $\ell$. Therefore, $\delta^+(H[X_{i-1}]) < c|X_{i-1}|$ by Lemma \ref{lemma:girth_3m}, meaning there is a vertex $u \in X_{i-1}$ with more than ${\delta^+(H) - c|X_{i-1}|}$ outneighbors outside of $X_{i-1}$, so
    \[|X_{i}| > \delta^+(H) + (1-c)|X_{i-1}|.\]
    This, by induction on $i$ and $|X_1| \geq \delta^+(H)$, gives 
    \[|X_i| \geq \delta^+(H)\cdot \frac{1 - (1-c)^i}{c} = \delta^+(H)\cdot \frac{1-0.6535^{i/m}}{1-0.6535^{1/m}}.
    \]
    In particular,
    $$
    |X_{5m}| \geq \delta^+(H) \cdot \frac{1-0.6535^{5}}{1-0.6535^{1/m}} = \frac nk\cdot \frac{1-0.6535^{5}}{1-0.6535^{1/m}} - o(n).
    $$
    We will show that for each $k$, $|X_{5m}| > \frac 12 n$.

    Since $m = \lfloor \frac{k-1}{3} \rfloor$, we have $k = 3m + r$ with $r \in \{1,2,3\}$. Consider the function 
    $$f(m) = \frac{1}{3m + r} \cdot \frac{1-0.6535^{5}}{1-0.6535^{1/m}}.$$ 
    By computing the derivative, one can verify that this function is increasing in $m$ for $m \geq 1$. Thus, in order to prove that $|X_{5m}| > \frac 12 n$ for each $k$, it is enough to check it for the smallest possible values of $k$ with different remainders modulo $3$. 
    For $k=4$ we have $|X_{5m}| \ge 0.635n - o(n)$, for $k=5$, $|X_{5m}| \ge 0.508n - o(n)$, while if $3 \mid k$, then $k\geq 9$ since $k$ must be a prime power, and $|X_{5m}| \ge 0.510n - o(n)$. In each case, $|X_{5m}| > \frac 12 n$ for large enough $n$.

    Therefore, for any vertex $v$, there is a directed path of length at most $5m$ from $v$ to over a half of the other vertices. This implies there has to be a vertex $v_0$ to which there are paths of length at most $5m$ from more than half of $V(H)$. 

    Since the set of vertices to which there is a path of length at most $5m$ from any vertex $v$ is larger than $\frac{1}{2}|V(H)|$, and the set of vertices from which there is a path of length at most $5m$ to $v_0$ is larger than $\frac{1}{2}|V(H)|$, their intersection is nonempty; concatenating the two paths through a common vertex we obtain a directed walk, and hence also a directed path, of length at most $10m$ from $v$ to $v_0$, which concludes the proof.
\end{proof}

Take the vertex $v_0 \in V(H)$ from Claim \ref{clm:short_paths_general} and define 
\begin{align*}
    &V_1 := N^+(v_0), \\
    &V_{i} := N^+(V_{i-1})\setminus \bigg(\bigcup_{j<i} V_j \bigg) \;\;\text{ for } i=2,3,\ldots,k.
\end{align*} 
In other words, for $i < k$, $V_i$ is the set of vertices at distance \textit{exactly} $i$ from $v_0$, while $V_k$ consists of the vertices at distance exactly $k$ from $v_0$ and possibly of $v_0$ itself (note that $v_0$ cannot lie in $V_j$ for $j<k$, as this would give a closed walk of length $j \leq k-1$). 

\begin{claim}\label{cla:homomorphic_to_c5}
There are no arcs from $V_i$ to $V_j$ in $H$ unless $j \equiv i+1 \bmod{k}$. 
\end{claim} 
\begin{proof}
Denote $D := 10\lfloor \frac{k-1}{3} \rfloor$.
If there was an arc $xy$ with $x\in V_i$, $y\in V_j$, then there would be paths of lengths $j$ and $i+1$ from $v_0$ to $y$. As there is a directed path of length $d \leq D$ from $y$ to $v_0$, we obtain intersecting closed directed walks of lengths $d+i+1$ and $d+j$ passing through $v_0$. Unless $j \equiv i+1 \bmod{k}$, we have $1 \leq |i+1-j| \leq k-1$, so $\gcd(d+i+1, d+j) \leq |i+1-j| < k$ and hence $\gcd(d+i+1, d+j) \mid \ell$, as every integer smaller than $k$ divides $\ell$. Since 
    $$\ell > 20k^2 > ((D+k+1)-1)((D+k)-1) \geq ((d+i+1)-1)(d+j-1),$$ 
$\ell$ is representable as $x(d+i+1) + y(d+j)$ with $x,y$ nonnegative integers (Corollary \ref{cor:frobenius}), and traversing the two closed walks $x$ and $y$ times, respectively, a forbidden member of $\WW{\ell}$ is obtained.
\end{proof}

As a result, we obtain $k$ sets, $V_1, \ldots, V_k$, each of which forms an independent set in $H$. Moreover, $|V_1| = d^+(v_0) \geq \delta^+(H)$ and for $i=2,3,\ldots,k$, $|V_i| \geq \delta^+(H) \geq \frac nk - o(n)$, since by Claim \ref{cla:homomorphic_to_c5} all outneighbors of a vertex of $V_{i-1}$ lie in $V_i$. Let $V' := V_1 \cup \ldots \cup V_k$. Then $|V'| \geq n - o(n)$. 

The graph $H[V']$ is homomorphic to $\CC{k}$ as needed, but it does not necessarily satisfy $\delta^{\pm}(H[V']) \geq \frac nk - o(n)$ yet. To achieve this, we have to remove some vertices from $V'$. 

\begin{claim}
    For any $\eps > 0$ there exists $\eps_1 > 0$ such that when $\delta^+(H) \geq \frac nk(1 - \eps_1)$, then one can remove at most $\eps n$ vertices from $V'$ so that the subgraph of $H$ induced by the remaining vertices has minimum semidegree at least $\frac nk(1 - \eps)$. In other words, one can remove $o(n)$ vertices from $V'$ so that the obtained induced subgraph has minimum semidegree at least $\frac nk - o(n)$. 
\end{claim}

\begin{proof}
Set $\eps_1 := \left(\frac{\eps}{4k}\right)^2$ and assume that $\delta^+(H) \geq \frac nk(1 - \eps_1)$.
Note that every vertex $v \in V(H)$ has inneighbors in at most one of the sets $V_i$. To see this, let $i$ be the smallest index such that $v$ has an inneighbor in $V_i$. If $i=k$, then $v$ has inneighbors only in $V_k$, while if $i \leq k-1$, then $v \in V_{i+1}$, and by Claim \ref{cla:homomorphic_to_c5} vertices from $V_{i+1}$ cannot have inneighbors in $V_j$ for $j \neq i$.

Since each $V_i$ is of size at least $\delta^+(H)$ and every $v\in V(H)$ can have inneighbors in at most one such set, 
\[
\Delta^-(H) \leq n - (k-1)\cdot \delta^+(H) \leq \delta^+(H) + \eps_1 n.
\]

Set $c := \frac{1}{\sqrt{\eps_1}}-1$ and partition $V(H)$ based on indegree in $H$:
\begin{align*}
V_{\text{low}} = \{ v \in V(H): d^-(v) \leq \delta^+(H) - c\eps_1 n \}, \\
V_{\text{high}} = \{ v \in V(H): d^-(v) > \delta^+(H) - c\eps_1 n \}.
\end{align*}
Then, 
\begin{align*}
|V(H)|\delta^+(H) &\leq \sum_{v\in V(H)} d^+(v) = \sum_{v\in V(H)} d^-(v) =
\sum_{v\in V_{\text{low}}} d^-(v) + \sum_{v\in V_{\text{high}}} d^-(v) \\
&\leq |V_{\text{low}}| \cdot (\delta^+(H) - c\eps_1 n) + |V_{\text{high}}| \cdot \Delta^-(H) \\ 
&\leq |V_{\text{low}}| \cdot (\delta^+(H) - c\eps_1 n) + |V_{\text{high}}| \cdot (\delta^+(H) + \eps_1 n) \\
&= |V(H)|\delta^+(H) + (|V_{\text{high}}| - c |V_{\text{low}}|)\cdot \eps_1 n \\
&\leq |V(H)|\delta^+(H) + (|V(H)| - (c+1) |V_{\text{low}}|)\cdot \eps_1 n.
\end{align*}
This implies that 
\[|V_{\text{low}}| \leq \frac{|V(H)|}{c+1} \leq n\sqrt{\eps_1}.\]

We can finally bound the minimum semidegree of $H[V' \setminus V_{\text{low}}]$ using the above bound on the size of $V_{\text{low}}$ and $|V(H) \setminus V'| \leq n - k\cdot \delta^+(H) \leq n\eps_1$:
\begin{align*}
\delta^+(H[V'\setminus V_{\text{low}}]) &\geq \delta^+(H) - |V(H) \setminus V'| - |V_{\text{low}}| \geq \delta^+(H) - n\eps_1 - n\sqrt{\eps_1} \\ 
&\geq \delta^+(H) - n \sqrt{\eps_1} - n \sqrt{\eps_1} \geq \frac nk (1 - \eps_1) - \frac{n\eps}{2k} \geq \frac nk (1 - \eps)
\end{align*}
and
\begin{align*}
\delta^-(H[V'\setminus V_{\text{low}}]) &\geq (\delta^+(H) - c\eps_1 n) - |V(H) \setminus V'| - |V_{\text{low}}| \geq \delta^+(H) - c\eps_1 n - n\eps_1 - n \sqrt{\eps_1} \\ 
&= \delta^+(H) - n\sqrt{\eps_1} - n\sqrt{\eps_1} \geq \frac nk (1 - \eps_1) - \frac{n\eps}{2k} \geq \frac nk(1 - \eps).
\end{align*}
Therefore, removing $V_{\text{low}}$ meets the requirement of the claim. 
\end{proof}

This concludes the proof of Lemma \ref{lemma:stability_even_case}.
\end{proof}

\subsection{Approximate structure in the odd case}\label{sec:odd_case_approx}

When $\ell$ is odd, the number $k$ from Theorem \ref{thm:main_thm} is equal to $4$. This leads to the following stability lemma, analogous to Lemma \ref{lemma:stability_even_case}.

\begin{lemma}\label{lemma:stability_odd_case}
    Let $\ell > 90$ be an odd integer divisible by $3$. Every $\CC{\ell}$-free oriented graph $G$ on $n$ vertices satisfying $\delta^+(G) > \frac n4$ contains a subgraph $H$ homomorphic to $\CC{4}$ with $\delta^{\pm}(H) \geq \frac n4 - o(n)$.
\end{lemma}

Applying Proposition \ref{prop:graph_cleanup} to the graph $G$, we obtain a $\WW{\ell}$-free subgraph $H$ with $\delta^+(H) \ge \frac n4 - o(n)$. However, in this case, $H$ does not necessarily contain a large induced subgraph homomorphic to $\CC{4}$ itself. Indeed, a random orientation of a sufficiently large balanced complete bipartite graph contains no odd cycles, has minimum outdegree $\frac n4 - o(n)$ and does not contain an induced subgraph on $n - o(n)$ vertices homomorphic to $\CC{4}$. This is why we need to prove the following weaker result first.

\begin{lemma}\label{lemma:stability_odd_case_1}
    Let $\ell > 90$ be an odd integer divisible by $3$. Every $\WW{\ell}$-free oriented graph $H$ on $n$ vertices with $\delta^+(H) \geq \frac n4 - o(n)$ contains an induced bipartite subgraph on $n-o(n)$ vertices with minimum semidegree at least $\frac n4 - o(n)$.
\end{lemma}
\begin{proof}

We start with a claim similar to Claim \ref{clm:short_paths_general}.

\begin{claim}\label{clm:bound_distance_odd_case}
    There exists a vertex $v_0 \in V(H)$ such that for every $v \in V(H)$, there is a directed path of length at most $6$ from $v$ to $v_0$.
\end{claim}
\begin{proof}
    Since $H$ is $\WW{\ell}$-free and $\ell$ is divisible by $3$, $H$ is $\CC{3}$-free. Take any $v \in V(H)$ and denote $X_1 = N^+(v)$, $X_{i+1} = N^+(X_i) \cup X_i$ for $i = 1,2$. Since for any $i$, $H[X_i]$ is $\CC{3}$-free, there must exist a vertex $u \in X_i$ such that $|N^+(u) \cap X_i| < 0.3465\cdot |X_i|$ (Lemma \ref{lemma:girth_3m} applied with $m=1$). Hence,
    $$
    |X_{i+1}| - |X_i| > \delta^+(H) - 0.3465 \cdot |X_i| = \frac n4 - o(n) - 0.3465\cdot |X_i|.
    $$
    This in particular implies that $|X_3| > 0.52n - o(n)$.
    In other words, every vertex $v$ has a directed path of length at most $3$ to over a half of $V(H)$. Therefore there must exist a vertex $v_0$ to which there is a directed path of length at most $3$ \emph{from} over a half of $V(H)$. Concatenating two such paths through a common vertex, we obtain a directed walk, and hence also a directed path, of length at most $6$ from any $v$ to $v_0$. 
\end{proof}

Take the vertex $v_0 \in V(H)$ from Claim \ref{clm:bound_distance_odd_case} and define
\begin{align*}
    &V_1 := N^+(v_0), \\
    &V_{i} := N^+(V_{i-1})\setminus \bigg(\bigcup_{j<i} V_j\bigg) \;\;\text{ for } i=2,3,4.
\end{align*} 
In other words, for $i \leq 3$, $V_i$ is the set of vertices at distance \emph{exactly} $i$ from $v_0$, while $V_4$ consists of the vertices at distance exactly $4$ and possibly of $v_0$ itself. We now observe that the following holds.

\begin{claim}\label{cla:odd_case_parity}
For any $i \equiv j \bmod 2$, there are no arcs from $V_i$ to $V_j$. 
\end{claim}
\begin{proof}
Assume there is an arc $xy$ with $x\in V_i$, $y\in V_j$, $i-j \in \{-2,0,2\}$. Then there are directed paths of lengths $i+1$ and $j$ from $v_0$ to $y$. We also know there is a directed path of length $d \leq 6$ from $y$ to $v_0$, which together with the previously mentioned paths forms intersecting closed directed walks of lengths $d+i+1$ and $d+j$. Now, since 
\[(d+i+1) - (d+j) = i-j+1 \in \{-1, 1, 3\},\] 
the number $\gcd(d+i+1, d+j)$ divides $3$, and hence it divides $\ell$ as well. Moreover, $\ell > 90 \geq (d+i)(d+j-1)$, therefore, according to Corollary \ref{cor:frobenius}, the closed directed walks of lengths $d+i+1$ and $d+j$ can be composed into a closed directed walk of length $\ell$, contradicting our assumptions on $H$.
\end{proof}

With careful analysis it is possible to reduce the bound to $\ell > 39$. Key observation is that there is also a cycle of length at most $6$ containing $v_0$. For the sake of brevity we omit the proof.

Claim \ref{cla:odd_case_parity} implies that the subgraph of $H$ induced by $V' := V_1 \cup V_2 \cup V_3 \cup V_4$ is bipartite with parts $V_1 \cup V_3$ and $V_2 \cup V_4$. We now prove that $|V'| \geq n - o(n)$.
Denote $T := V(H) \setminus V'$. 

\begin{claim}\label{clm:odd_case_small_trash_set}
For any $\eps > 0$ there exists $\eps_1 > 0$ such that $\delta^+(H) \geq \frac n4(1 - \eps_1)$ implies $|T| \leq \eps n$. In particular, $|T| = o(n)$. 
\end{claim}
\begin{proof}
Take $\eps_1 < \frac{3\eps^2}{1 + 3\eps}$ and assume that $\delta^+(H) \geq \frac n4 (1 - \eps_1)$. First bound the number of arcs leaving $V'$ as below.
$$
\delta^+(H)\cdot (n-|T|)\leq\sum_{v\in V'}d^+(v)= e(V',V')+e(V',T) \leq \frac{(n-|T|)^2}{4}+|V_4||T|.
$$
The last inequality is due to $H[V']$ being bipartite and arcs leaving $V'$ solely from $V_4$, as a vertex with an inneighbor in $V_i$ for some $i \leq 3$ is at distance at most $4$ from $v_0$ and thus belongs to $V'$. This can be rewritten as
$$
\frac n4 (1-\eps_1) \leq \delta^+(H) \leq \frac{n-|T|}{4}+|V_4|\cdot \frac{|T|}{n-|T|}.
$$

Assume now that $|T| > \eps n$, aiming for a contradiction. We bound $|V_4|$ using the inequality above:
\begin{align*}
|V_4| &\geq \frac{|T|-\eps_1n}{4}\cdot \frac{n-|T|}{|T|} = \frac{n-|T|}{4} - \frac{n-|T|}{4|T|}\cdot \eps_1n > \frac{n-|T|}{4} - \frac{\eps_1}{4\eps}n \\
&= \frac 14(|V_1|+|V_2|+|V_3|+|V_4|) - \frac{\eps_1}{4\eps}n,
\end{align*}
which implies that
\[|V_4| > \frac 13(|V_1|+|V_2|+|V_3|) - \frac{\eps_1}{3\eps}n.\]

We now focus on bounding $|V_1|+|V_2|+|V_3|$. Notice that all arcs from $V_1$ leave to $V_2$ and all arcs from $V_2$ leave either back to $V_1$ or to $V_3$. Moreover, as $H$ is an oriented graph, between a vertex of $V_1$ and a vertex of $V_2$ there is at most one arc, hence
$$ \delta^+(H)\cdot (|V_1|+|V_2|) \leq e(V_1,V_2) + e(V_2,V_1\cup V_3) \leq |V_1||V_2| + |V_2||V_3| = |V_2|(|V_1|+|V_3|) $$
and therefore, using $|V_1| \geq \delta^+(H)$,
\[|V_1|+|V_3| \geq \delta^+(H) \cdot \frac{|V_1|}{|V_2|}+\delta^+(H) \geq \frac{(\delta^+(H))^2}{|V_2|} +\delta^+(H).\] 
Thus, by the inequality of arithmetic and geometric means,
\[|V_1|+|V_2|+|V_3| \geq \frac{(\delta^+(H))^2}{|V_2|} +\delta^+(H) + |V_2| \geq 3 \delta^+(H)\]
and
\[|V'| > 4\delta^+(H) - \frac{\eps_1}{3\eps}n \geq 4 \cdot \frac n4(1 - \eps_1) - \frac{\eps_1}{3\eps}n = n\left(1 - \eps_1 \cdot \frac{3\eps + 1}{3\eps}\right) > n(1 - \eps).\]
This is a contradiction, as by our previous assumptions, $|V'| = n - |T| \leq n - \eps n$.
\end{proof}

The graph $H[V']$ is large and bipartite as wanted, but it does not need to satisfy $\delta^\pm(H[V']) \geq \frac n4 - o(n)$ yet. We need to remove some vertices from $V'$.

\begin{claim}\label{clm:stability_outdegree_odd_case}
    Let $\eps > 0$. There exists $\eps_1 > 0$ such that if $\delta^+(H) \geq \frac n4(1 - \eps_1)$, then one can remove at most $\eps n$ vertices from $V'$ so that the subgraph of $H$ induced by the remaining vertices has minimum semidegree at least $\frac n4 - 3 \eps n$. In other words, one can remove $o(n)$ vertices from $V'$ so that the obtained induced subgraph has minimum semidegree at least $\frac n4 - o(n)$.
\end{claim}
\begin{proof}
    Denote $n' := |V'|= n - |T|$. By Claim \ref{clm:odd_case_small_trash_set}, there exists $\eps_2 > 0$, for which $\delta^+(H) \geq \frac n4 (1 - \eps_2)$ implies that $|T| \leq \frac{\eps^2}{8}n$, so $n' \geq n - \frac{\eps^2}{8} n$. 
    
    Take $\eps_1 := \min (\eps_2, \frac{\eps^2}{2})$. Then $\delta^+(H) \geq \frac n4(1 - \eps_1)$ implies 
    \begin{align*}
    &n' \geq n - \frac{\eps^2}{8} n \geq n (1 - \eps), \\
    &\delta^+(H[V']) \geq \frac n4 (1 - \eps_1) - |T| \geq \frac n4 (1 - \eps^2) \geq \frac{n'}{4}(1 - \eps^2).
    \end{align*}
    
    Consider now a set $Y$ of $y = \eps n'$ vertices of the smallest indegree in $H[V']$. Denote $s := \max_{v \in Y} \{ d^-_{H[V']}(v) \}$. Then $\delta^-(H[V'\setminus Y]) \geq s - y$. 
    Calculating the number of arcs of $H[V']$ from $V' \setminus Y$ we obtain the following.
    \begin{align*}
    (n'-y)\cdot \frac{n'}{4}(1 -\eps^2) &\leq \sum_{v\in V'\setminus Y} d_{H[V']}^+(v) = e(V'\setminus Y, V'\setminus Y) + e(V'\setminus Y, Y) \\ &\leq \frac{(n'-y)^2}{4} + y s. 
    \end{align*}  
Here the first inequality follows from the bound on $\delta^+(H[V'])$ above, while the last one is a consequence of the bipartiteness of $H[V'\setminus Y]$ together with the definition of $s$. This allows us to lower-bound $s$.
$$
    s \geq \frac{(n'-y)(y - \eps^2 n')}{4y} = \frac{(n'-\eps n')(\eps n' - \eps^2 n')}{4\eps n'} = \frac{n'(1-\eps)^2}{4} > \frac{n'}{4}(1-2\eps).
$$
This means that $\delta^-(H[V' \setminus Y]) \geq s - y > \frac{n'}4 - 2\eps n'$. Clearly, after removing $\eps n'$ vertices, $\delta^+(H[V' \setminus Y]) \geq \delta^+(H[V']) - \eps n' > \frac{n'}{4} - 2 \eps n'$, therefore 
$$
\delta^\pm(H[V'\setminus Y]) \geq n'\left(\frac{1}{4} - 2\eps\right) \geq n(1 - \eps)\left(\frac 14 - 2\eps\right) \geq \frac n4 - 3 \eps n
$$
as needed.
\end{proof}

By Claim \ref{clm:stability_outdegree_odd_case}, after removing set $Y$ of vertices from $V'$, the graph $H[V' \setminus Y]$ satisfies the conditions of Lemma \ref{lemma:stability_odd_case_1}.
\end{proof}

Before we proceed to the proof of Lemma~\ref{lemma:stability_odd_case}, we make a general observation that will be useful later.

\begin{obs}\label{obs:odd_case_anchor}
Let $G$ be an oriented graph on $n$ vertices with $\delta^+(G) > \frac n4$, and let $V(G) = L \cup R \cup T$ be a partition with $|T| = o(n)$. For $n$ large enough, there exists a directed path $P$ in $G$ of length at most $3$ with endpoints $u,v$ in $L \cup R$ and all internal vertices in $T$, such that the length of $P$ is odd if $u$ and $v$ lie in the same part of the partition $(L,R)$ and even otherwise.
\end{obs}

\begin{proof}
Assume that no such path exists. Then, there are no arcs inside $L$ or $R$, because any such arc would be a wanted path satisfying all the requirements. 
Partition $T$ into $T_L \cup T_R \cup T_0$, where
\begin{align*}
& T_L = \{t\in T: N^-(t)\cap R \not= \emptyset\}, \\
& T_R = \{t\in T: N^-(t)\cap L \not= \emptyset\}, \\
& T_0 = T\setminus (T_L\cup T_R).
\end{align*}
Notice that $T_L$ and $T_R$ are disjoint, because a vertex $t$ lying in both would have inneighbors in both $L$ and $R$ and, by the outdegree assumption, an outneighbor in $L \cup R$, creating a wanted path of length $2$. For the same reason, for every $t \in T_L$, we have $N^+(t) \cap L = \emptyset$, so $|N^+(t) \cap R| \geq \frac n4 - o(n)$. Consequently, there are no arcs in $T_L$, because an arc $tt'$ inside $T_L$, together with $u \in N^-(t)\cap R$ and $v \in N^+(t')\cap R$, would give a path $utt'v$ of length $3$ with both endpoints in $R$. Symmetrically for $T_R$. 

Consider the partition of $V(G) \setminus T_0$ into $L \cup T_L$ and $R \cup T_R$. By the above, every arc of $G - T_0$ joins the two parts.
Moreover no arc leads from $L \cup R$ to $T_0$, so at most $|T||T_0|$ arcs leave $G - T_0$. Hence, by summing the outdegrees of the vertices of $G-T_0$ we obtain,
\[\frac{n}{4}\left(n-|T_0|\right) < \sum_{v \in V(G) \setminus T_0} d^+(v) \leq e(G - T_0) + |T||T_0| \leq \frac{(n-|T_0|)^2}{4} + |T||T_0|.\]
This leads to $\frac{1}{4}|T_0|(n-|T_0|-4|T|) < 0$, which is false for large $n$, as $|T_0| \leq |T| = o(n)$.
\end{proof}

\begin{proof}[Proof of Lemma~\ref{lemma:stability_odd_case}]
Take the graph $G$ and use Proposition \ref{prop:graph_cleanup} and Lemma \ref{lemma:stability_odd_case_1} to obtain a bipartite subgraph $H$ with $\delta^\pm(H) \ge \frac n4 - o(n)$. This allows us to partition $V(G)$ into $L$, $R$ (the parts of $H$), and the remaining vertex set $T$ (short for `trash'). Moreover, $|T| = o(n)$. Note there may be extra arcs in $G[L \cup R]$ which violate the bipartition $(L, R)$. 

The path $P$ from Observation \ref{obs:odd_case_anchor}, applied to the partition $(L, R, T)$, will serve as our anchor. Assume without loss of generality that $u \in L$. Greedily extend $P$ beyond $v$ inside $H$, avoiding the vertices chosen before (the condition $\delta^\pm(H) \geq \frac n4 - o(n)$ guarantees this is always possible), to obtain a directed path $Q$ of length $\ell - 3$ from $u$ to some vertex $w$. Every arc of $H$ joins $L$ with $R$, so the part containing $w$ is obtained from the part containing $v$ by $\ell - 3 - |P|$ swaps, while by the choice of $P$ the part containing $v$ is obtained from the one containing $u$ by $|P| - 1$ swaps. Altogether $w$ lies in the part obtained from $L$ by $\ell - 4$ swaps, that is, $w \in R$. Denote
\begin{align*}
&V_1 = N^+_{H}(w) \setminus Q,  &V_4 = N^-_{H}(u) \setminus Q.
\end{align*}

Note that $V_1 \subseteq L$ and $V_4 \subseteq R$. Since $\delta^\pm(H) \geq \frac n4 - o(n)$, $|V_1| \geq \frac n4 - o(n)$ and $|V_4| \geq \frac n4 - o(n)$.
Observe moreover that there are no arcs from $V_1$ to $V_4$ in $G$, as such an arc would close the path $Q$ into a $\CC{\ell}$. Now let
\begin{align*}
&V_2 = R \setminus V_4,  &V_3 = L \setminus V_1 .
\end{align*}
Since there are no arcs from $V_1$ to $V_4$ and $\delta^\pm(H) \geq \frac n4 - o(n)$, we infer that $V_2 \supseteq N^+_{H}(V_1)$ has at least $\frac n4 - o(n)$ vertices and similarly that $|V_3| \geq \frac n4 - o(n)$. This also implies that for any $i$, $|V_i| \leq \frac n4 + o(n)$. In other words, the partition $V_1, V_2, V_3, V_4$ is nearly balanced. Observe that
\begin{align*}
    &e(V_1, V_2) = e(V_1, R) - e(V_1, V_4) = e(V_1, R) \geq |V_1| \cdot \delta^+(H) \geq \frac{n^2}{16} - o(n^2), \\
    &e(V_2, V_1) \leq |V_2||V_1| - e(V_1, V_2) \leq \left(\frac n4 + o(n)\right)^2 - \left(\frac{n^2}{16} - o(n^2)\right) \leq o(n^2), \\
    &e(V_2, V_3) = e(V_2, L) - e(V_2, V_1) \geq |V_2| \cdot \delta^+(H) - o(n^2) \geq \frac{n^2}{16} -o(n^2). 
\end{align*}

This reasoning can be continued until we obtain:
\begin{align*}
    &e(V_1, V_2) = \frac{n^2}{16} - o(n^2)  &e(V_2, V_1) = o(n^2) \\
    &e(V_2, V_3) = \frac{n^2}{16} - o(n^2)  &e(V_3, V_2) = o(n^2) \\
    &e(V_3, V_4) = \frac{n^2}{16} - o(n^2)  &e(V_4, V_3) = o(n^2) \\
    &e(V_4, V_1) = \frac{n^2}{16} - o(n^2)  &e(V_1, V_4) = o(n^2) 
\end{align*}

We found a partition of $V':=V_1 \cup V_2 \cup V_3 \cup V_4$ which is close to the blob structure of a blow-up of $\CC{4}$, but we do not know anything about its minimum semidegree yet. We need to once again refine the partition.

Denote $b := e_{H}(V_2, V_1) + e_{H}(V_3, V_2) + e_{H}(V_4, V_3) + e_{H}(V_1, V_4)$ (the number of `backward arcs' in $H$). We know that $b = o(n^2)$.

\begin{claim}\label{clm:odd_case_final_cleanup}
    Let $\eps > 0$. There exists $\eps_1 > 0$ such that whenever $\delta^\pm(H) \geq \frac n4(1 - \eps_1)$ and $b \leq \eps_1 n^2$, one can remove at most $\eps n$ vertices and at most $\eps n^2$ arcs from $H$ so that the obtained graph $H'$ is homomorphic to $\CC{4}$ and $\delta^\pm(H') \geq \frac n4 - \eps n$. 
    In other words, one can remove $o(n^2)$ arcs and $o(n)$ vertices from $H$ so that the obtained graph is homomorphic to $\CC{4}$ and has minimum semidegree at least $\frac n4 - o(n)$. 
\end{claim}
\begin{proof}
    Let $\eps_1 = \frac{\eps^2}{16}$. We analyze the arcs which do not agree with the $V_1, \ldots, V_4$ blob structure. Take any vertex $v \in V(H)$ and let $v \in V_i$. Let $b_1(v)$ denote the number of arcs from $v$ to $V_{i-1}$ and $b_2(v)$ denote the number of arcs from $V_{i+1}$ to $v$. Since every backward arc is counted once in $\sum_{v \in V'} b_1(v)$ and once in $\sum_{v \in V'} b_2(v)$, we know that $\sum_{v\in V'} (b_1(v) + b_2(v)) = 2b \leq 2\eps_1 n^2$. 

    Consider the set $T'$ of $\frac{\eps n}{2}$ vertices with the largest values of $b(v) := b_1(v) + b_2(v)$. Denoting $s := \min_{v \in T'} b(v)$, we have
    $$
    \frac{\eps n}{2} \cdot s \leq \sum_{v \in T'} b(v) \leq \sum_{v\in V'} b(v) \leq 2\eps_1 n^2 = \frac{\eps^2}{8}n^2,
    $$
    It follows that $s \leq \frac{\eps n}{4}$.

    We now consider the subgraph $H'$ on vertex set $V' \setminus T'$ and with all arcs of kind $V_i \rightarrow V_{i+1}$. Clearly, $H'$ is homomorphic to $\CC{4}$. Moreover, every vertex $v \in V' \setminus T'$ satisfies $b(v) \leq s$, which gives
    $$
    \delta^\pm(H') \geq \delta^\pm(H) - |T'| - s \geq \frac n4 - \frac{\eps_1 n}{4} - \frac{\eps n}{2} - \frac{\eps n}{4} \geq \frac n4 - \eps n.
    $$

    We removed at most $\frac{\eps}{2} n < \eps n$ vertices and at most $b \leq \eps n^2$ arcs from $H$, thus conditions of the claim are met.
\end{proof}

The graph $H'$ from Claim~\ref{clm:odd_case_final_cleanup} satisfies the conditions of Lemma \ref{lemma:stability_odd_case} and finishes its proof. 
\end{proof}

\subsection{Framework for building cycles}\label{section:framework}

For any $\ell > 20k^2$ divisible by $3$, consider a $\CC{\ell}$-free oriented graph $G$ with $\delta^+(G) > \frac nk$. By Lemma~\ref{lemma:stability_even_case} or Lemma~\ref{lemma:stability_odd_case} we know that $G$ contains a (not necessarily induced) subgraph $H$ homomorphic to $\CC{k}$ satisfying $\delta^\pm(H) \geq \frac nk - o(n)$. Denote the blobs of $H$ as $V_1, \ldots, V_k$ in the cyclic order and $V' := V_1 \cup \ldots \cup V_k$. Since each blob contains the outneighborhood in $H$ of a vertex of the preceding blob, we have $|V_i| \geq \frac nk - o(n)$ for every $i$, and hence also $|V_i| \leq n - (k-1)\left(\frac nk - o(n)\right) = \frac nk + o(n)$.

Consider how $H$ is embedded in $G$. The graph $G[V']$ might contain arcs which do not go from $V_i$ to $V_{i+1}$. Call these \textit{extra arcs}. Moreover, let $T := V(G) \setminus V'$ denote the set of \textit{extra vertices}. As $|V'| \geq n - o(n)$, we know that $|T| = o(n)$.

Unless specified otherwise, neighborhoods are considered in the graph $G$.

\begin{obs}\label{obs:large_intersection}
For any $x \in V_i$, $y \in V_{i+2}$ (indices taken modulo $k$), 
$$
|N^+(x) \cap N^-(y) \cap V_{i+1}| \geq \frac nk - o(n).
$$
\end{obs}
\begin{proof}
    Since $N_{H}^+(x), N_{H}^-(y) \subseteq V_{i+1}$, $d_{H}^+(x), d_{H}^-(y) \geq \frac nk - o(n)$, and $|V_{i+1}| \leq \frac nk + o(n)$, 
    \begin{align*}
    |N_{H}^+(x) \cap N_{H}^-(y) \cap V_{i+1}| &\geq d_{H}^+(x) + d_{H}^-(y) - |V_{i+1}| \\ &\geq \left(\frac nk - o(n)\right) + \left(\frac nk - o(n)\right) - \left(\frac nk + o(n)\right) = \frac nk - o(n).
    \end{align*}
    Since $H$ is a subgraph of $G$, the same bound holds for the neighborhoods taken in $G$.
\end{proof}

\begin{obs}\label{obs:closing_short_paths}
For any $x \in V_i$, $y\in V_j$, $x\neq y$, there exists a directed path from $x$ to $y$ in $H$ of length at most $k+1$. 
\end{obs}
\begin{proof} 
Take $x$ and greedily build a sequence $x_0 = x, x_1 \in N^+_{H}(x_0), x_2 \in N^+_{H}(x_1), \ldots$, until we reach a vertex $x'$ in $V_{j-2}$ (clearly, this will happen after at most $k-1$ steps). Now apply Observation \ref{obs:large_intersection} to conclude there is a path of length $2$ from $x'$ to $y$, which, together with the previous vertices, forms a walk of length at most $k+1$. Observe that at each step we had at least $\frac nk - o(n)$ vertices to choose from, which means we can easily ensure that no vertex repeats on the directed walk.
\end{proof}

The degree condition $\delta^+(G) > \frac nk$ forces $G$ to contain elements which do not follow the blow-up structure of $H$. These can be used when traversing the blobs of $H$ to build a cycle of length not divisible by $k$. The definition below formalizes this concept. 

\begin{defn}\label{def:detour}
    A \textit{detour} is any directed path $P$ in $G$ which starts in $V_i$ and ends in $V_j$ for some indices $i, j$. 
    If the length of $P$ is $d$, then the \textit{remainder} of $P$ is $i-j+d \bmod k$. A detour is \textit{non-zero} if its remainder is not $0$.
\end{defn}

Note that an extra arc is also a non-zero detour. The following observation will allow us to connect detours in order to obtain a cycle of length $\ell$.

\begin{obs}\label{obs:joining_detours}
Let $\mathcal{P}$ be a collection of less than $k$ pairwise vertex-disjoint detours in $G$, whose sum of lengths is at most $3k$. If the sum of remainders of elements of $\mathcal{P}$ is $\ell \bmod k$, then the elements of $\mathcal{P}$ can be combined with arcs of $H$ into a $\CC{\ell}$.
\end{obs}
\begin{proof}
    Let $\mathcal{P} = \{P_1, P_2, \ldots, P_a\}$. Denote the first vertex on the path $P_i$ as $u_i$ and the last vertex of $P_i$ as $v_i$. Now apply Observation \ref{obs:closing_short_paths} to the pairs $(v_1, u_2), (v_2, u_3), \ldots, (v_a, u_1)$ to connect the paths in $\mathcal{P}$. Recall that for each of the `linking' vertices we had at least $\frac nk - o(n)$ possibilities to choose from. In particular, we can pick them to be pairwise distinct and disjoint from $\mathcal{P}$. We obtain a directed cycle $C$, denote its length by $d$. Since each linking path lies in $H$, its length is congruent modulo $k$ to the difference of the indices of the blobs it joins, and therefore $d \equiv \ell \bmod k$. Moreover, 
    \[d \leq |\mathcal{P}|\cdot (k+1) + 3k < k^2 + 4k < \ell,\]
    because for even $\ell$ we have $\ell > 20k^2$, while for odd $\ell$ we know that $k=4$ and $\ell > 90 > 32$. It now suffices to extend $C$ until it has length $\ell$. 
    
    Perform as many times as needed the following operation. Find in $C$ an arc $xy$ which agrees with the blob structure (i.e.~belongs to $H$); such an arc must exist due to how we constructed the linking paths in Observation \ref{obs:closing_short_paths}. Now replace the arc $xy$ in $C$ with the sequence $xw_1w_2\ldots w_ky$, where $w_1 \in N_{H}^+(x) \setminus C, w_2 \in N^+_{H}(w_1) \setminus C$, etc. This clearly extends the length of the cycle $C$ by $k$.
\end{proof}

This observation allows us to get rid of all extra arcs, simplifying the structure of $G[V']$. Note that an extra arc $xy$ with $x\in V_i, y\in V_j$ constitutes a non-zero detour of remainder $i-j+1 \bmod k$.

\begin{obs}\label{obs:reducing_extra_edges}
    By moving at most $2(k-2)$ vertices from $V'$ to $T$, one can ensure that every arc of $G[V']$ leads from some $V_i$ to $V_{i+1}$. The resulting partition still satisfies $\delta^{\pm}(G[V']) \geq \frac nk - o(n)$ and $|T| = o(n)$.
\end{obs}

\begin{proof}
Assume first that $G$ contains $k-1$ pairwise vertex-disjoint extra arcs. Each of them creates a detour of non-zero remainder. 
If $\ell$ is even, then by Corollary \ref{cor:remainders_plain} we can choose a subset of these arcs whose corresponding remainders sum up to $\ell \bmod k$. These are fewer than $k$ pairwise vertex-disjoint detours, so by Observation \ref{obs:joining_detours} the graph $G$ would contain a $\CC{\ell}$, a contradiction. 
Hence, consider odd $\ell$, so $k=4$, and check how the remainders of the $3$ extra arcs can be distributed. None of the remainders is equal to $\ell \bmod 4$, as otherwise we have a contradiction by Observation~\ref{obs:joining_detours}. If all three remainders are equal to $\ell + 2 \bmod 4$, then their sum is $\ell \bmod 4$ and again we have a contradiction by Observation~\ref{obs:joining_detours}. Hence, there is a remainder $2$. If any of the remaining remainders is $\ell + 2 \bmod 4$, then the two sum up to $\ell \bmod 4$, contradicting Observation~\ref{obs:joining_detours}. Thus, all three extra arcs are creating detours of remainder $2$. Apply now Observation~\ref{obs:odd_case_anchor} to the partition of $V(G)$ into $V_1 \cup V_3$, $V_2 \cup V_4$ and $T$. The path $P$ it provides is a detour, and its remainder is odd. This remainder is not $\ell \bmod 4$, as otherwise Observation~\ref{obs:joining_detours} applied to $P$ alone gives a $\CC{\ell}$, so it equals $\ell + 2 \bmod 4$. As the internal vertices of $P$ lie in $T$, only the two endpoints of $P$ can meet the six endpoints of our three extra arcs, so at least one of these arcs is vertex-disjoint from $P$. The two remainders sum up to $\ell \bmod 4$, contradicting Observation~\ref{obs:joining_detours} as previously. 

Therefore a maximum collection of pairwise vertex-disjoint extra arcs consists of at most $k-2$ arcs, and by maximality every extra arc of $G$ meets one of its at most $2(k-2)$ endpoints. Moving these endpoints from $V'$ to $T$ thus destroys all extra arcs. As only a constant number of vertices is moved, the bounds on $|T|$ and on the semidegrees of $G[V']$ are unaffected.
\end{proof}

From now on we assume that every arc of $G[V']$ leads from some $V_i$ to $V_{i+1}$. Recall that the set of extra vertices $T = V(G) \setminus V'$ satisfies $|T| = o(n)$. If a vertex $v \in T$ has an arc from some $V_i$ and an arc to $V_j$ ($i, j$ not necessarily distinct), then it forms a detour of length $2$ with remainder $i-j+2 \bmod k$. We introduce the following definition to be able to easily connect these detours.

\begin{defn}\label{def:in_pointing}
    A blob $V_i$ is an \emph{in-pointing} blob for a vertex $t \in T$ if $|N^-(t) \cap V_i| \geq 2k$. Similarly, a blob $V_i$ is \emph{out-pointing} for some $t\in T$ if $|N^+(t) \cap V_i| \geq 2k$.
\end{defn}

This neighborhood size bound will prove useful when expanding an obtained closed directed walk of length $\ell$ into a proper $\CC{\ell}$.

Every $t \in T$ has at most $2k^2$ arcs into blobs which are not out-pointing for $t$, and similarly, $t$ has at most $2k^2$ arcs from blobs which are not in-pointing for $t$. Summing this up, we obtain at most $4k^2|T| = o(n)$ vertices in $V'$ which can be adjacent to $T$, but not within the out- or in-pointing blob relation. Call this set of exceptional vertices $X$. Note that for any $i$, $|V_i \setminus X| \geq \frac nk - o(n)$. 

\begin{defn}\label{def:sidewalk}
A \emph{sidewalk} $V_i \rightarrow t \rightarrow V_j$ is a tuple $(V_i, t, V_j)$ with $t \in T$, in which $V_i$ is in-pointing for $t$ and $V_j$ is out-pointing for $t$. The vertex $t$ is the \textit{center} of this sidewalk. The \emph{remainder} of a sidewalk $V_i \rightarrow t \rightarrow V_j$ is $i-j+2 \bmod k$. A sidewalk is \emph{non-zero} if its remainder is not $0$.
\end{defn}

Note that each sidewalk corresponds to many potential detours. The observation below allows us to join sidewalks to a cycle of length $\ell$. In this section we used it only for joining sidewalks, but in Section~\ref{section:special_cases} it will be used also to join sidewalks to a single longer detour.

\begin{obs}\label{obs:join_sidewalks}
Let $\mathcal{S}$ be a collection of sidewalks and at most one detour, such that in total there are fewer than $k$ elements, centers of all sidewalks are distinct and different from the vertices of the detour, and the detour has length at most $k$ and at most $4$ vertices from $V'$. If the remainders of the elements of $\mathcal{S}$ add up to $\ell \bmod k$, then $G$ contains a $\CC{\ell}$. 
\end{obs}
\begin{proof}
    For each sidewalk $V_i \rightarrow t \rightarrow V_j$ in $\mathcal{S}$, choose some vertices $u_t \in V_i \cap N^-(t)$ and $w_t \in V_j \cap N^+(t)$, so that all the chosen vertices $u_t, w_t$ are pairwise distinct and different from the detour. This is possible since for each sidewalk, $|V_i \cap N^-(t)| \geq 2k$ and $|V_j \cap N^+(t)| \geq 2k$, while at each step fewer than $2k$ vertices have to be avoided. As a result, we obtain less than $k$ vertex disjoint detours of total length at most $2(k-2)+k < 3k$, whose remainders add up to $\ell \bmod k$, which allows us to conclude using Observation \ref{obs:joining_detours}.
\end{proof}

The proven observation allows us to bound the total number of non-zero sidewalks. 
Note that a vertex $t \in T$ may be the center of several sidewalks with the same remainder, but we want to count every remainder only once for every center. 

\begin{obs}\label{obs:sidewalks_bound}
    For every $t \in T$, let $S_t$ be the set of non-zero remainders of the sidewalks with center $t$. Then $\sum_{t \in T} |S_t| \leq k-2$.
\end{obs}

\begin{proof}
Assume there is a set of $a$ centers $t_1, t_2, \ldots, t_a \in T$ and $a$ non-empty sets of remainders $R_1, \ldots, R_a \subseteq (\mathbb{Z}_k \setminus \{0\})$ such that $t_i$ forms non-zero sidewalks of all the remainders in $R_i$, and $\sum_{i=1}^a |R_i| \geq k-1$. Discarding remainders if necessary, we may assume that $\sum_{i=1}^a|R_i| = k-1$ and hence $a \leq k-1$. 

For even $\ell$, by Lemma \ref{lemma:remainders_sets} we are able to choose at most one remainder from each $R_i$ so that the chosen remainders sum up to $\ell \bmod k$. The corresponding sidewalks have distinct centers and there are fewer than $k$ of them. The obtained set of sidewalks can be used to build a $\CC{\ell}$ by Observation~\ref{obs:join_sidewalks}, a contradiction. 

For odd $\ell$, so $k=4$, repeatedly using Observation \ref{obs:join_sidewalks} we deduce how the remainders are distributed. Firstly, no sidewalk has remainder $\ell \bmod 4$, so all sidewalks are of remainders $2$ or $\ell+2 \bmod 4$ and consequently $a \geq 2$. Then, there cannot be sidewalks with distinct centers that have remainders $2$ and $\ell+2 \bmod 4$, so all sidewalks have the same remainder and $a = 3$. Finally, there are no three sidewalks of remainder $\ell+2 \bmod 4$, thus there are three sidewalks of remainder $2$ with distinct centers. Let $P$ be the detour given by Observation~\ref{obs:odd_case_anchor} applied to the partition of $V(G)$ into $V_1 \cup V_3$, $V_2 \cup V_4$ and $T$. As in the proof of Observation~\ref{obs:reducing_extra_edges} its remainder is $\ell + 2 \bmod 4$. Since $P$ has at most two internal vertices, the center of one of the previously obtained three sidewalks, say $t \in T$, does not lie on $P$. By the definition of a sidewalk, $t$ has at least $8$ inneighbors and outneighbors in the respective blobs, so it creates a detour that is vertex-disjoint from $P$. The remainders of this detour and $P$ sum to $\ell \bmod 4$, which by Observation \ref{obs:joining_detours} yields a $\CC{\ell}$, a contradiction.    
\end{proof}

%

\subsection{Conclusion}\label{sec:main_proof_finish}

We can now  finish the proof of Theorem~\ref{thm:main_thm}.

\begin{proof}[Proof of Theorem~\ref{thm:main_thm}]
Assuming the contrary, let $G$ be a large enough oriented graph on $n$ vertices with $\delta^+(G) \geq \frac nk + \frac{k-1}{k}$ not containing a directed cycle of length $\ell > 20k^2$ divisible by $3$. 
By Lemma~\ref{lemma:stability_even_case} or Lemma~\ref{lemma:stability_odd_case} and Observation~\ref{obs:reducing_extra_edges}, the vertex set of $G$ can be partitioned into sets $V' = V_1 \cup \ldots \cup V_k$ and $T$, such that every arc of $G[V']$ leads from some $V_i$ to $V_{i+1}$,  $\delta^\pm(G[V']) \geq \frac nk - o(n)$ and $|T|=o(n)$.
We will show that the sets $S_t$ of non-zero remainders of the sidewalks with center $t$ are too large, which leads to a contradiction with Observation~\ref{obs:sidewalks_bound}.

For each $t \in T$, denote by $I_t$ the set of indices $i$ such that $V_i$ is in-pointing for $t$. Every $t \in T$ has at least one out-pointing blob: indeed $d^+(t) \geq \delta^+(G) - |T| \geq \frac nk - o(n)$, so at least one of the $k$ blobs receives at least $\frac{n}{k^2} - o(n) \geq 2k$ arcs from $t$. Fixing such a blob $V_j$, the $|I_t|$ remainders $i-j+2$ of the sidewalks $V_i \rightarrow t \rightarrow V_j$ are pairwise distinct and at most one of them equals $0$. Hence $|S_t| \geq |I_t|-1$ for every $t \in T$.

Consider a sample of $k$ vertices $w_1 \in V_1 \setminus X, w_2 \in V_2 \setminus X, \ldots, w_k \in V_k \setminus X$. For any $t\in T$, it holds that if there is an arc $w_it$, then $V_i$ is in-pointing for $t$. Moreover, since every arc of $G[V']$ leads from some $V_i$ to $V_{i+1}$, we know that $N^+(w_i) \subseteq (V_{i+1} \cup T)$. Therefore,
\begin{align}\label{eq:counting_sidewalks_main}
n + k - 1 &\leq k \cdot \delta^+(G) \leq \sum_{i=1}^k d^+(w_i)
= \sum_{i=1}^k \left( |N^+(w_i) \cap V_{i+1}| + |N^+(w_i) \cap T| \right) \nonumber \\
&\leq \sum_{i=1}^k |V_{i+1}| + \sum_{t\in T} |I_t|
= (n - |T|) + \sum_{t \in T} |I_t| = n + \sum_{t \in T} (|I_t|-1).
\end{align}

This implies that $\sum_{t \in T} |S_t| \geq \sum_{t \in T} (|I_t|-1) \geq k - 1$. This contradicts Observation~\ref{obs:sidewalks_bound} and finishes the proof of Theorem \ref{thm:main_thm}.
\end{proof}


\section{Special cases} \label{section:special_cases}

Theorem~\ref{thm:main_thm} is not tight for four families of lengths, for which the optimal outdegree bound is determined by Theorem~\ref{thm:special_cases}. Together with the constructions of Section~\ref{section:constructions}, Theorem~\ref{thm:special_cases} is directly implied by the following theorem. 

\begin{theorem}\label{thm:exceptional}
    Let $\ell$ be an integer divisible by $3$ and $k$ be the smallest integer greater than $2$ that does not divide $\ell$, such that $(k,\ell \bmod k) \in \{(4,2),\, (4,3),\, (5,3),\, (5,4)\}$. Assume that $\ell > 20k^2$ and put
    \[c := \begin{cases} k-3 & \text{if } \ell \bmod k = 3, \\ k-2 & \text{if } \ell \bmod k \neq 3 .\end{cases}\]
    Then every oriented graph $G$ on $n$ vertices, with $n$ large enough, satisfying $\delta^+(G) \geq \frac{n+c}{k}$ contains a directed cycle of length $\ell$.
\end{theorem}

Throughout the section we assume, aiming at a contradiction, that $G$ is a $\CC{\ell}$-free oriented graph on $n$ vertices with $\delta^+(G) \geq \frac{n+c}{k}$. In particular $\delta^+(G) > \frac nk$, so we may use everything proved in Subsection~\ref{section:framework}. We fix the partition $V(G) = V' \cup T$ with ${V' = V_1 \cup \ldots \cup V_k}$ provided by Lemma~\ref{lemma:stability_even_case} or Lemma~\ref{lemma:stability_odd_case} together with Observation~\ref{obs:reducing_extra_edges}, so that every arc of $G[V']$ leads from some $V_i$ to $V_{i+1}$, $\delta^\pm(G[V']) \geq \frac nk - o(n)$ and $|T| = o(n)$. We keep the in-pointing and out-pointing relations of Definition~\ref{def:in_pointing}, the sidewalks of Definition~\ref{def:sidewalk} and the exceptional set $X \subseteq V'$ with $|X| = o(n)$.  Recall that any adjacency between a vertex $t \in T$ and a vertex of $V_i \setminus X$ witnesses the corresponding pointing relation, and that $|V_i \setminus X| \geq \frac nk - o(n)$. 

\subsection{Admissible sidewalk remainders}\label{ss:exc_count}

Contrary to the proof of Theorem~\ref{thm:main_thm}, now the inequalities in~(\ref{eq:counting_sidewalks_main}) with the weaker bound $\delta^+(G) \geq \frac{n+c}{k}$ will not provide $k-1$ non-zero sidewalks. We need to do more careful analysis. Because of that, we define quantities that measure how ``irregular'' the vertices in $T$ and the blobs in $V'$ behave.

For $t \in T$ we define $I_t$ as the set of indices $i$ such that $V_i$ is in-pointing for $t$, and $O_t$ as the set of indices $j$ such that $t$ has an outneighbor in $V_j \setminus X$. Note that this is a bit stronger property than just $V_j$ is an out-pointing blob for $t$.
Nevertheless, by the definition of $X$, every such $V_j$ is out-pointing for $t$. 
Let $R_t$ be the set of non-zero remainders of the sidewalks $V_i \rightarrow t \rightarrow V_j$ with $i \in I_t$ and $j \in O_t$, that is 
\[R_t= (I_t - O_t + 2)\setminus \{0\},\]
and let
\[\eps_t = |R_t| - |I_t| + 1.\] 
Since $t$ has at least $\delta^+(G) - |T| - |X| > 0$ outneighbors in $V' \setminus X$, the set $O_t$ is non-empty, so $\eps_t \ge 0$ and $\eps_t$ measures how many additional non-zero remainders $t$ creates. 
In particular, $\eps_t=0$ implies that $I_t - j + 2 = R_t \cup \{0\}$ for every $j \in O_t$, that is, $t$ together with $V_j$ and all its in-pointing blobs creates sidewalks of all remainders in $R_t \cup \{0\}$.

While for $i \in \{1,2,\ldots,k\}$ let $T_i$ be the set of vertices in $T$, for which $V_i$ is an in-pointing blob, that is,
\[T_i = \{t \in T : i \in I_t\}\]
and let 
\[\mu_i = |V_{i+1}| + |T_i| - \delta^+(G).\] 
Since each vertex in $V_i \setminus X$ has all its outneighbors in $V_{i+1}$ or in $T_i$, we have $\mu_i \geq 0$, and $\mu_i$ measures how much room is left. In particular, $\mu_i = 0$ implies that every vertex in $V_i \setminus X$ has an arc to every vertex in $V_{i+1}$ and to all vertices $t \in T$, for which $V_i$ is in-pointing. 

Summing up the defined measures over all blobs and vertices in $T$, we obtain 
\begin{align*}
0 &\leq \sum_{i=1}^k \mu_i + \sum_{t\in T} \eps_t = \Big( n - |T| + \sum_{t\in T} |I_t| - k\delta^+(G) \Big) + \Big( \sum_{t\in T} |R_t| - \sum_{t\in T} |I_t|  + |T|\Big) \\
&= \sum_{t\in T} |R_t| + n - k\delta^+(G) \leq \sum_{t\in T} |R_t| - c.
\end{align*}
Thus, $\sum_{t\in T} |R_t| \geq c$. On the other hand, as $R_t \subseteq S_t$, by Observation~\ref{obs:sidewalks_bound}, we have $\sum_{t \in T} |R_t| \leq k-2$. Since $c=k-2$ or $k-3$, we have that 
\begin{equation}\label{eq:totaldefect}
     \sum_{i=1}^k \mu_i + \sum_{t\in T} \eps_t \in \{0,1\}
\end{equation}
and $\sum_{t\in T} |R_t|$ is equal to $k-2$ or $k-3$.
Note that $\sum_{t\in T}|R_t| = k-3$ is possible only when $c = k-3$, that is when $\ell \bmod k = 3$, and it forces all $\mu_i=0$ and all $\eps_t = 0$. 
Therefore, depending on the considered value of $\ell$ and $k$, and the value of $\sum_{t\in T}|R_t|$, we have $6$ possible cases. They are listed in Table~\ref{tab:special_cases}.
Before we analyze those cases in detail, we make some observations on possible sidewalk remainders that allow us to list all possible non-empty sets $R_t$ in each case. 

\begin{obs}\label{obs:2sidewalk}
If $i \in I_t \cap O_t$ for some $t \in T$, that is, $2 \in R_t$, then $\mu_i = 1$.
\end{obs}

\begin{proof}
By the definition of $O_t$, there exists an arc $tv$ for some vertex $v \in V_i \setminus X$. Since all outneighbors of $v$ lie in $V_{i+1}$ or in $T_i \setminus\{t\}$, 
\[\delta^+(G) \leq d^+(v) \leq |V_{i+1}| +  |T_i| - 1 = \delta^+(G) + \mu_i - 1,\]
so $\mu_i \geq 1$ and $\mu_i = 1$ by \eqref{eq:totaldefect}. 
\end{proof}

\begin{obs}\label{obs:eps0}
If $\eps_t = 0$ for some $t \in T$, then $|O_t|=1$. 
\end{obs}

\begin{proof}
Assumption $\eps_t = 0$ means that for every $j \in O_t$ there are $|I_t| = |R_t|+1$ distinct remainders $i-j+2$ in $R_t \cup \{0\}$, so $I_t - j + 2 = R_t \cup \{0\}$. This implies that if there are two distinct $j, j' \in O_t$, then $R_t \cup \{0\} = (R_t \cup \{0\}) + (j'-j)$, so $R_t \cup \{0\}$ contains the subgroup generated by $j'-j$. Among the considered cases, there is only one non-trivial subgroup of $\mathbb{Z}_k$ that does not contain $\ell \bmod k$, that is $\{0,2\}$ in the case $k=4$ and $\ell \bmod k = 3$. But then, $R_t \cup \{0\}$ contains $\{0,2\}$, is contained  in $\{0,1,2\}$, and is invariant under the translation by $2$, so it is equal to $\{0,2\}$. Thus, $R_t=\{2\}$, $O_t = \{j, j+2\}$, and the equality $I_t - j + 2 = R_t \cup \{0\}$ gives $I_t = \{j, j+2\}$. By Observation~\ref{obs:2sidewalk}, $\mu_j = 1$ and $\mu_{j+2} = 1$ contradicting \eqref{eq:totaldefect}. 
\end{proof}

The above observations exclude the possibility of having two distinct vertices $t, t' \in T$ with $2 \in R_t \cap R_{t'}$. Indeed, in such a case $i \in I_t \cap O_t$ and $i' \in I_{t'} \cap O_{t'}$ for some $i, i'$, so by Observation~\ref{obs:2sidewalk}, $\mu_i=\mu_{i'}=1$, and \eqref{eq:totaldefect} implies that $i=i'$ and $\eps_t=\eps_{t'}=0$. Thus, by Observation~\ref{obs:eps0}, $O_t = O_{t'} = \{i\}$. Since the outneighbors of $t$ and $t'$ are only in $T$, $X$ and $V_i \setminus X$, each of them has at least $\delta^+(G) - |T| - |X| = \frac nk - o(n)$ outneighbors in $V_i \setminus X$, while $|V_i| \leq \frac nk + o(n)$. Hence, there is a vertex $v \in V_i \setminus X$ with both arcs $tv$ and $t'v$. All outneighbors of $v$ lie in $V_{i+1}$ or in $T_i \setminus\{t,t'\}$, so
\[\delta^+(G) \leq d^+(v) \leq |V_{i+1}| +  |T_i|-2 = \delta^+(G) + \mu_i -2,\]
which means that $\mu_i \geq 2$, contradicting \eqref{eq:totaldefect}. 

Therefore, there is at most one vertex $t \in T$ with $2 \in R_t$, and it may exist only when $\sum_i \mu_i = 1$, which means $c=k-3$ and $\sum_t |R_t|=k-2$. 
This, together with Observation~\ref{obs:join_sidewalks}, allows us to list all possible sets of remainders in each considered case. 

\begin{table}[h!]
    \begin{tabular}{ccc@{\hspace{10mm}}l}
        $(k,\ell \bmod k)$ & $\sum_i \mu_i + \sum_t \eps_t$ & $\sum_t |R_t|$ & possible families of non-zero remainders \\[3pt] \hline
        $(4,2)$ & 0 & $k-2=2$ & $(\{1,3\})$,$\quad (\{1\},\{3\})$ \\[3pt]
        $(4,3)$ & 0 & $k-3=1$ & $(\{1\})$ \\[3pt]
        $(4,3)$ & 0 or 1 & $k-2=2$ & $(\{1\},\{1\})$,$\quad (\{1,2\})$ \\[3pt]
        $(5,3)$ & 0 & $k-3=2$ & $(\{1\},\{1\})$, $\quad (\{1\},\{4\})$, $\quad (\{1,4\})$ \\[3pt]
        \multirow{2}{*}{$(5,3)$} & \multirow{2}{*}{0 or 1} & \multirow{2}{*}{$k-2=3$} & $(\{1\},\{1\},\{4\})$, $\quad (\{1\},\{1,4\})$,  \\
         & & & $(\{1,2\},\{4\})$, $\quad (\{1,2,4\})$ \\[3pt]
        $(5,4)$ & 0 & $k-2=3$ & $(\{1\},\{1\},\{1\})$ 
    \end{tabular}
    \caption{All possible families of the non-empty sets $R_t$.}\label{tab:special_cases}
\end{table}

Determination of the above possible families is straightforward. As an example, we show the most complex case $(k,\ell \bmod k) = (5,3)$ with $\sum_t|R_t| = k-2 = 3$. 
Assume first that there is a vertex $t \in T$ having remainder $2$ in $R_t$. As shown above, no other vertex has remainder $2$. By Observation~\ref{obs:join_sidewalks}, no other vertex in $T$ has remainder $1$, no vertex has remainder $3$, and there are no two vertices having remainder $4$, so in order to have $3$ non-zero remainders, we need to have a single vertex with $R_t = \{1,2,4\}$ or one vertex with $\{1,2\}$ and the other with $\{4\}$. 
Whereas, if there is no vertex with remainder $2$, then we may have at most two vertices with remainder $1$ and at most one with remainder $4$, which leads to two possible families depending if the vertices having remainders $1$ and $4$ are distinct. 

By \eqref{eq:totaldefect} and Observation~\ref{obs:eps0}, all but at most one vertex $t \in T$ satisfy $|O_t|=1$. This allows us to ``assign'' every such vertex to an appropriate blob in order to have most of its outneighbors in the consecutive blob.
Formally, for any vertex $t \in T$ with $O_t = \{i\}$, we define the \emph{position} of $t$ as $p(t)=i-1$. 
We also say that any vertex $v \in V_i$ has position $p(v)=i$. Thus, all but at most one vertex of $G$ (possibly the vertex $t \in T$ having $\eps_t=1$) has a position.
Note that every vertex with a position $i$ has all its outneighbors in $V_{i+1} \cup X \cup T$. 
If $\eps_t = 0$ for some $t \in T$, then $I_t - p(t) + 1 = R_t \cup \{0\}$, so $V_{p(t)-1}$ is in-pointing for $t$, and $t \in T_{p(t)}$ if and only if $1 \in R_t$.

We will build detours from short paths of vertices with positions. For a directed path $v_1 v_2 \ldots v_m$ with $m \leq 3$ of vertices in $T \cup X$ with positions, let $V_i$ be an in-pointing blob of $v_1$ if $v_1 \in T$, or $i = p(v_1)-1$ if $v_1 \in X$. Adding an inneighbor of $v_1$ in $V_i$ at the beginning and, if $v_m \in T$, an outneighbor of $v_m$ in $V_{p(v_m)+1}$ at the end, we obtain a detour of remainder $i - p(v_m) + m$ with at most four vertices in $V'$, and each added vertex can be chosen in at least $2k$ ways. We say that $v_1, \ldots, v_m$ \emph{create a detour} of this remainder.

\subsection{Finalization in the cases with $\sum_i \mu_i + \sum_t \eps_t = 0$}

The following observation will be crucial and used multiple times. 

\begin{obs}\label{obs:special_cases_escape}
If $\sum_i \mu_i + \sum_t \eps_t = 0$ and $v$ is a vertex such that at least $x$ vertices in $V_{p(v)+1} \cup T_{p(v)}$ are not outneighbors of $v$, then $v$ has at least $x$ outneighbors outside $V_{p(v)+1}  \cup T_{p(v)}$, each with a position different than $p(v)+1$. In particular, every vertex $v \in T$ satisfying $1 \in R_v$ has an outneighbor outside $V_{p(v)+1}  \cup T_{p(v)}$ with a position different than $p(v)+1$.
\end{obs}

\begin{proof}
Since $\mu_{p(v)}=0$,
\[d^+(v) \geq \delta^+(G) = |V_{p(v)+1}| + |T_{p(v)}|,\]
and so $v$ has at least $x$ outneighbors outside $V_{p(v)+1} \cup T_{p(v)}$.
Every such outneighbor in $V'$ has a position different than $p(v)+1$, as wanted. If there is an outneighbor $u$ in $T \setminus T_{p(v)}$ with position $p(v)+1$, then, since $\eps_u = 0$, $V_{p(v)}$ is an in-pointing blob of $u$, contradicting the fact that $u \not\in T_{p(v)}$. 

If $v$ satisfies $1 \in R_v$, then $v$ is in $T_{p(v)}$, so not all vertices in $T_{p(v)}$ are its outneighbors and $v$ has an outneighbor with a position different than $p(v)+1$.
\end{proof}
 
This allows us to solve all the cases with $\sum_i \mu_i + \sum_t \eps_t = 0$.
Each of the cases below is illustrated in a figure. There, the blobs are labelled relative to the position of $v$, the thick gray arrows represent the arcs between consecutive blobs. An arrow from a blob to a vertex, or from a vertex to a blob, indicates that the vertex has at least $2k$ inneighbors, respectively outneighbors, in this blob. All the drawn arrows are forced by the arguments, and the highlighted paths are the detours and sidewalks leading to a contradiction.

\smallskip
\noindent\textbf{Case $(k,\ell \bmod k) = (4,2)$:}
\smallskip

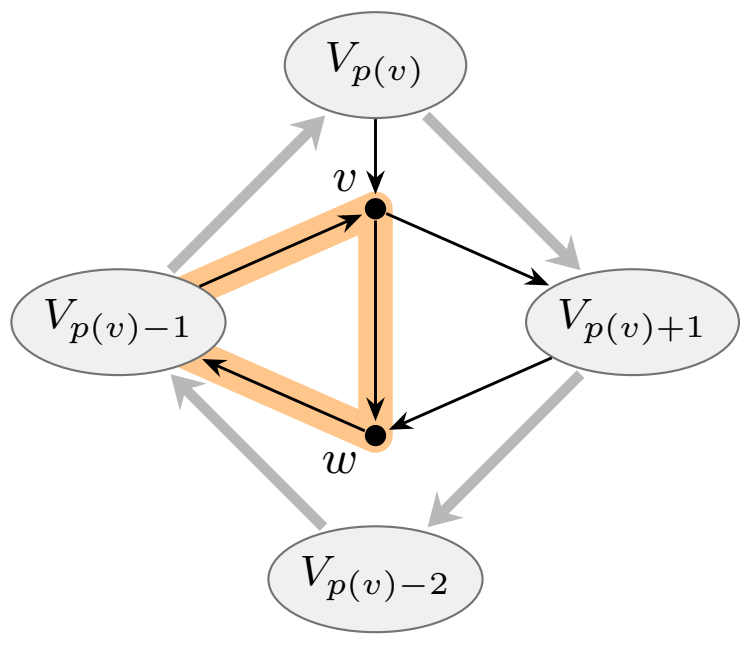
\begin{figure}[ht]
\centering
\begin{tikzpicture}
  \blobsfour
  \node[vtx,label={[vlab]above left:$v$}] (v) at (0,0.75) {};
  \node[vtx,label={[vlab]below left:$w$}] (w) at (0,-0.75) {};
  \begin{scope}[on background layer]
    \draw[hl] (Bm1.center) -- (v.center) -- (w.center) -- (Bm1.center);
  \end{scope}
  \draw[arc] (Bm1) -- (v);
  \draw[arc] (B0) -- (v);
  \draw[arc] (v) -- (B1);
  \draw[arc] (B1) -- (w);
  \draw[arc] (w) -- (Bm1);
  \draw[arc] (v) -- (w);
\end{tikzpicture}
\caption{Case $(k,\ell \bmod k)=(4,2)$.}\label{fig:case42}
\end{figure}

Let $v \in T$ be a vertex with $1 \in R_v$. By Observation~\ref{obs:special_cases_escape}, $v$ has an outneighbor $w$ outside $T_{p(v)}$ with a position different than $p(v)+1$. To avoid having a detour of remainder $2 \bmod 4$, the position of $w$ must be also different than $p(v)$ and $p(v)-1$, so $p(w) = p(v)-2$. This means that $v$ and $w$ creates a detour of remainder $3 \bmod 4$ (see Figure~\ref{fig:case42}). To avoid contradiction with Observation~\ref{obs:join_sidewalks}, the sidewalk of remainder $3$ must have a center $v$ or $w$. In the former, $V_{p(v)-2}$ is an in-pointing blob of $v$, which creates with $w$ a forbidden detour of remainder $2 \bmod 4$. While in the latter, $V_{p(v)}$ is an in-pointing blob of $w$, so $w \in T_{p(v)}$, contradicting the definition of $w$.

\smallskip
\noindent\textbf{Case $(k,\ell \bmod k) = (4,3)$ with $\sum_i \mu_i + \sum_t \eps_t = 0$:}
\smallskip

\begin{figure}[ht]
\centering
\begin{tikzpicture}
  \blobsfour
  \node[vtx,label={[vlab]above left:$v$}] (v) at (0,0.95) {};
  \node[vtx,label={[vlab]above left:$w$}] (w) at (0,0.35) {};
  \node[vtx,label={[vlab]below left:$u$}] (u) at (0,-0.3) {};
  \begin{scope}[on background layer]
    \draw[hl] (B0.center) -- (v.center) -- (w.center) -- (u.center) -- (B1.center);
  \end{scope}
  \draw[arc] (Bm1) -- (v);
  \draw[arc] (B0) -- (v);
  \draw[arc] (v) -- (B1);
  \draw[arc] (Bm1) -- (w);
  \draw[arc] (w) -- (B1);
  \draw[arc] (Bm1) -- (u);
  \draw[arc] (u) -- (B1);
  \draw[arc] (v) -- (w);
  \draw[arc] (w) -- (u);
\end{tikzpicture}
\hfill
\begin{tikzpicture}
  \blobsfour
  \node[vtx,label={[vlab]above right:$v$}] (v) at (0,0.95) {};
  \node[vtx,label={[vlab]below:$w$}] (w) at (0.5,0.3) {};
  \node[vtx,label={[vlab]below left:$u$}] (u) at (-0.75,-0.45) {};
  \begin{scope}[on background layer]
    \draw[hl] (Bm1.center) -- (v.center) -- (w.center) -- (u.center) -- (B0.center);
  \end{scope}
  \draw[arc] (Bm1) -- (v);
  \draw[arc] (B0) -- (v);
  \draw[arc] (v) -- (B1);
  \draw[arc] (Bm1) -- (w);
  \draw[arc] (w) -- (B1);
  \draw[arc] (Bm2) -- (u);
  \draw[arc] (u) -- (B0);
  \draw[arc] (v) -- (w);
  \draw[arc] (w) -- (u);
\end{tikzpicture}
\hfill
\begin{tikzpicture}
  \blobsfour
  \node[vtx,label={[vlab]above left:$v$}] (v) at (0,0.95) {};
  \node[vtx,label={[vlab]above left:$w$}] (w) at (0,0.35) {};
  \node[vtx,label={[vlab]below left:$u$}] (u) at (0,-0.75) {};
  \begin{scope}[on background layer]
    \draw[hl] (Bm1.center) -- (w.center) -- (u.center) -- (Bm1.center);
  \end{scope}
  \draw[arc] (Bm1) -- (v);
  \draw[arc] (B0) -- (v);
  \draw[arc] (v) -- (B1);
  \draw[arc] (Bm1) -- (w);
  \draw[arc] (w) -- (B1);
  \draw[arc] (B1) -- (u);
  \draw[arc] (u) -- (Bm1);
  \draw[arc] (v) -- (w);
  \draw[arc] (w) -- (u);
\end{tikzpicture}
\caption{Case $(k,\ell \bmod k)=(4,3)$ with $\sum_i \mu_i + \sum_t \eps_t = 0$ and all possible positions of~$u$.}\label{fig:case43}
\end{figure}
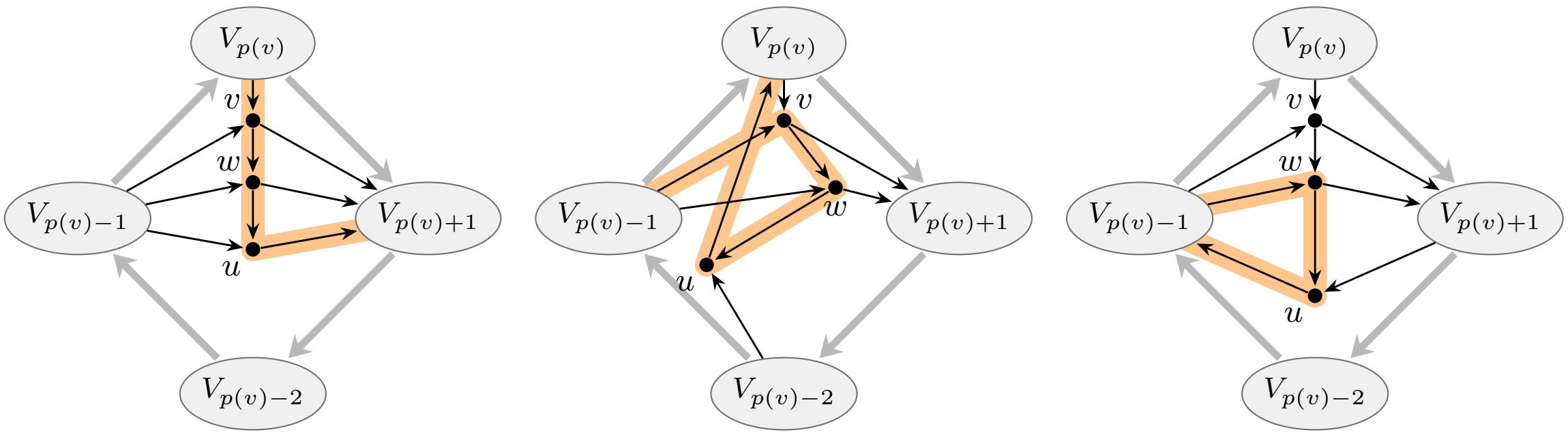

Similarly as in the previous case, let $v \in T$ be a vertex with $1 \in R_v$. By Observation~\ref{obs:special_cases_escape}, $v$ has an outneighbor $w$ outside $T_{p(v)}$ with a position different than $p(v)+1$. To avoid having a detour of remainder $3 \bmod 4$, the position of $w$ must be also different than $p(v)-1$ and $p(v)-2$, so $p(w) = p(v)$. 
Applying again Observation~\ref{obs:special_cases_escape}, but for vertex $w$, for which $v \in T_{p(v)}$ is not an outneighbor, we obtain an outneighbor $u$ of $w$ outside $T_{p(v)}$ with a position different than $p(v)+1$. Independently of the position of $u$, there is a detour of remainder $3 \bmod 4$ (see Figure~\ref{fig:case43}). 

\smallskip
\noindent\textbf{Case $(k,\ell \bmod k) = (5,3)$ with $\sum_i \mu_i + \sum_t \eps_t = 0$:}
\smallskip

Let $v \in T$ be a vertex with $1 \in R_v$. By Observation~\ref{obs:special_cases_escape}, $v$ has an outneighbor $w$ outside $T_{p(v)}$ with a position different than $p(v)+1$. To avoid having a detour of remainder $3 \bmod 5$, the position of $w$ must be also different than $p(v)-1$ and $p(v)-2$, so $p(w)=p(v)$ or $p(w)=p(v)+2$. We consider those two options separately. Also distinguishing between the case of family $(\{1\}, \{1\})$ and all the other families. Note that other families always include a sidewalk of remainder $4$. 

Consider first all the families except $(\{1\}, \{1\})$. If $p(w)=p(v)$, then we apply Observation~\ref{obs:special_cases_escape} for vertex $w$, for which $v$ is not an outneighbor, and obtain a vertex $u$ outside $T_{p(v)}$ with a position different than $p(v)+1$. The only option avoiding a detour of remainder $3 \bmod 5$ is to have $p(u)=p(v)+2$. Then $w$ and $u$ create a detour of remainder $4 \bmod  5$ (see Figure~\ref{fig:case53}(a)). By Observation~\ref{obs:join_sidewalks}, we can combine it with a sidewalk of remainder $4$, unless $4 \in R_w$ or $4 \in R_u$. In the former $w$ and $u$ create a detour of remainder $3 \bmod 5$, while in the latter $u \in T_{p(v)}$ contradicting its definition. 
If $p(w)=p(v)+2$, then $v$ and $w$ create a detour of remainder $4 \bmod  5$ (see Figure~\ref{fig:case53}(b)). We can combine it with a sidewalk of remainder $4 \bmod 5$, unless $4 \in R_v$ or $4 \in R_w$. Similarly as earlier, in the former $v$ and $w$ create a detour of remainder $3 \bmod 5$, while in the latter $w \in T_{p(v)}$ contradicting its definition. 

Therefore, consider now the family $(\{1\}, \{1\})$. If $p(w)=p(v)$, then $v$ and $w$ create a detour of remainder $2$. We can combine it with the second center of a sidewalk of remainder~$1$, unless $1 \in R_w$. But then $p(v)$ is an in-pointing blob of $w$, so $w \in T_{p(v)}$ contradicting its definition. Thus, $p(w)=p(v)+2$ and $v$ and $w$ create a detour of remainder $4 \bmod  5$. 
We can use the same reasoning for the second vertex $v'$ satisfying $1 \in R_{v'}$ to obtain vertex $w'$ with $p(w')=p(v')+2$ creating a detour of remainder $4 \bmod  5$. If $\{v,w\} \cap \{v',w'\} = \emptyset$, then, choosing the added vertices distinct, we obtain two vertex-disjoint detours of remainder $4$, which can be combined using Observation~\ref{obs:joining_detours} into a $\CC{\ell}$.
If $v'=w$, then $v$, $w$ and $w'$ create a detour of remainder $3 \bmod 5$ (see Figure~\ref{fig:case53}(c)), and symmetrically if $w' = v$. While if $w'=w$, then $p(v)=p(v')$, and one of the arcs $vv'$ and $v'v$ is missing, say $v'$ is not an outneighbor of $v$. Since $v,v' \in T_{p(v)}$, by Observation~\ref{obs:special_cases_escape}, $v$ has an outneighbor $u$ distinct from $w$ outside $V_{p(v)+1} \cup T_{p(v)}$ with a position different than $p(v)+1$. As for $w$, we get $p(u)=p(v)+2$, so $v$ and $u$ create a detour of remainder $4$. Since $u \neq v'$ and $w \neq v$, the pairs $\{v,u\}$ and $\{v',w\}$ are disjoint, which leads to a contradiction as before (see Figure~\ref{fig:case53}(d)). 

\begin{figure}[ht]
\centering
\begin{tikzpicture}
  \blobsfive
  \node[vtx,label={[vlab]above left:$v$}] (v) at (0,1.05) {};
  \node[vtx,label={[vlab]below left:$w$}] (w) at (0,0.4) {};
  \node[vtx,label={[vlab]right:$u$}] (u) at (0.56,-0.77) {};
  \begin{scope}[on background layer]
    \draw[hl] (Bm1.center) -- (w.center) -- (u.center) -- (Bm2.center);
  \end{scope}
  \draw[arc] (Bm1) -- (v);
  \draw[arc] (B0) -- (v);
  \draw[arc] (v) -- (B1);
  \draw[arc] (Bm1) -- (w);
  \draw[arc] (w) -- (B1);
  \draw[arc] (B1) -- (u);
  \draw[arc] (u) -- (Bm2);
  \draw[arc] (v) -- (w);
  \draw[arc] (w) -- (u);
  \node[pcap] at (0,-2.35) {(a) $p(w)=p(v)$};
\end{tikzpicture}
\hspace{6mm}
\begin{tikzpicture}
  \blobsfive
  \node[vtx,label={[vlab]above left:$v$}] (v) at (0,0.85) {};
  \node[vtx,label={[vlab]right:$w$}] (w) at (0.5,-0.69) {};
  \begin{scope}[on background layer]
    \draw[hl] (Bm1.center) -- (v.center) -- (w.center) -- (Bm2.center);
  \end{scope}
  \draw[arc] (Bm1) -- (v);
  \draw[arc] (B0) -- (v);
  \draw[arc] (v) -- (B1);
  \draw[arc] (B1) -- (w);
  \draw[arc] (w) -- (Bm2);
  \draw[arc] (v) -- (w);
  \node[pcap] at (0,-2.35) {(b) $p(w)=p(v)+2$};
\end{tikzpicture}
\\[4mm]
\begin{tikzpicture}
  \blobsfive
  \node[vtx,label={[vlab]above right:$v$}] (v) at (0,0.85) {};
  \node[vtx,label={[vlab]right:$w$}] (w) at (0.5,-0.69) {};
  \node[vtx,label={[vlab]left:$w'$}] (wp) at (-0.95,-0.25) {};
  \begin{scope}[on background layer]
    \draw[hl] (Bm1.center) -- (v.center) -- (w.center) -- (wp.center) -- (B0.center);
  \end{scope}
  \draw[arc] (Bm1) -- (v);
  \draw[arc] (B0) -- (v);
  \draw[arc] (v) -- (B1);
  \draw[arc] (B1) -- (w);
  \draw[arc] (B2) -- (w);
  \draw[arc] (w) -- (Bm2);
  \draw[arc] (Bm2) -- (wp);
  \draw[arc] (wp) -- (B0);
  \draw[arc] (v) -- (w);
  \draw[arc] (w) -- (wp);
  \node[pcap] at (0,-2.35) {(c) $w=v'$};
\end{tikzpicture}
\hspace{6mm}
\begin{tikzpicture}
  \blobsfive
  \node[vtx,label={[vlab]above right:$v$}] (v) at (0,1.05) {};
  \node[vtx,label={[vlab]below left:$v'$}] (vp) at (-0.45,0.35) {};
  \node[vtx,label={[vlab]below:$w$}] (w) at (0.35,-0.35) {};
  \node[vtx,label={[vlab]right:$u$}] (u) at (0.85,-0.65) {};
  \begin{scope}[on background layer]
    \draw[hl] (Bm1.center) -- (v.center) -- (u.center) -- (Bm2.center);
    \draw[hlb] (Bm1.center) -- (vp.center) -- (w.center) -- (Bm2.center);
  \end{scope}
  \draw[arc] (Bm1) -- (v);
  \draw[arc] (B0) -- (v);
  \draw[arc] (v) -- (B1);
  \draw[arc] (Bm1) -- (vp);
  \draw[arc] (B0) -- (vp);
  \draw[arc] (vp) -- (B1);
  \draw[arc] (B1) -- (w);
  \draw[arc] (w) -- (Bm2);
  \draw[arc] (B1) -- (u);
  \draw[arc] (u) -- (Bm2);
  \draw[arc] (v) -- (u);
  \draw[arc] (v) -- (w);
  \draw[arc] (vp) -- (w);
  \node[pcap] at (0,-2.35) {(d) $w=w'$};
\end{tikzpicture}
\caption{Case $(k,\ell \bmod k)=(5,3)$ with $\sum_i \mu_i + \sum_t \eps_t = 0$: (a), (b) for the families containing a sidewalk of remainder~$4$; (c), (d) for the family $(\{1\},\{1\})$ with $p(w)=p(v)+2$ and $p(w')=p(v')+2$.}\label{fig:case53}
\end{figure}

\smallskip
\noindent\textbf{Case $(k,\ell \bmod k) = (5,4)$:}
\smallskip

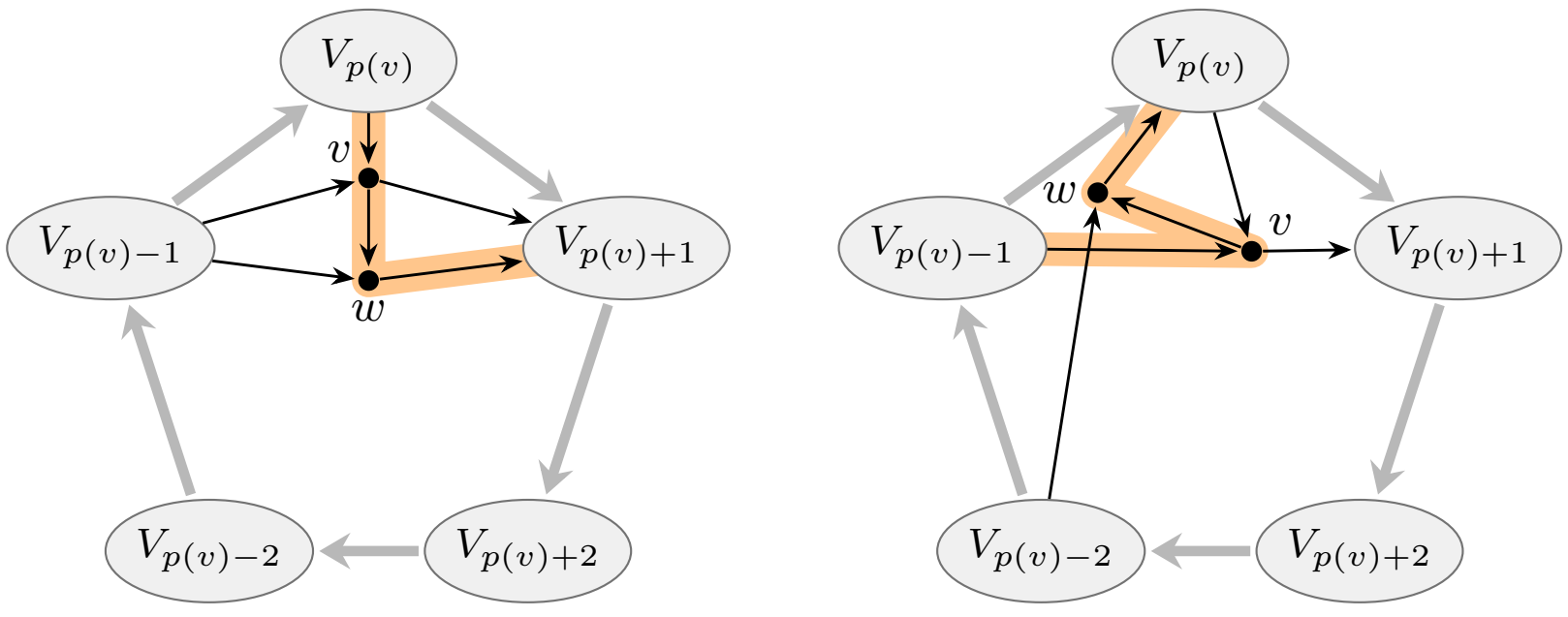
\begin{figure}[h!t]
\centering
\begin{tikzpicture}
  \blobsfive
  \node[vtx,label={[vlab]above left:$v$}] (v) at (0,1.05) {};
  \node[vtx,label={[vlab]below:$w$}] (w) at (0,0.35) {};
  \begin{scope}[on background layer]
    \draw[hl] (B0.center) -- (v.center) -- (w.center) -- (B1.center);
  \end{scope}
  \draw[arc] (Bm1) -- (v);
  \draw[arc] (B0) -- (v);
  \draw[arc] (v) -- (B1);
  \draw[arc] (Bm1) -- (w);
  \draw[arc] (w) -- (B1);
  \draw[arc] (v) -- (w);
  \node[pcap] at (0,-2.35) {(a) $p(w)=p(v)$};
\end{tikzpicture}
\hspace{6mm}
\begin{tikzpicture}
  \blobsfive
  \node[vtx,label={[vlab]above right:$v$}] (v) at (0.35,0.55) {};
  \node[vtx,label={[vlab]left:$w$}] (w) at (-0.7,0.95) {};
  \begin{scope}[on background layer]
    \draw[hl] (Bm1.center) -- (v.center) -- (w.center) -- (B0.center);
  \end{scope}
  \draw[arc] (Bm1) -- (v);
  \draw[arc] (B0) -- (v);
  \draw[arc] (v) -- (B1);
  \draw[arc] (Bm2) -- (w);
  \draw[arc] (w) -- (B0);
  \draw[arc] (v) -- (w);
  \node[pcap] at (0,-2.35) {(b) $p(w)=p(v)-1$};
\end{tikzpicture}
\caption{Case $(k,\ell \bmod k)=(5,4)$ with the two possible positions of $w$.}\label{fig:case54}
\end{figure}

Among all the three vertices $t \in T$ with $1 \in R_t$, let $v$ be such that $p(v)-1$ is different than the position of the other two vertices. 
Similarly as in the previous cases, by Observation~\ref{obs:special_cases_escape}, $v$ has an outneighbor $w$ outside $T_{p(v)}$ with a position different than $p(v)+1$. To avoid having a detour of remainder $4 \bmod 5$, the position of $w$ must be also different than $p(v)-2$ and $p(v)+2$, so $p(w) = p(v)-1$ or $p(w) = p(v)$ (see Figure~\ref{fig:case54}). 
If $p(w) = p(v)$, then $v$ and $w$ create a detour of remainder $2$ that can be combined by Observation~\ref{obs:join_sidewalks} with the sidewalks of remainder $1$ using other two vertices in $T$, unless $1 \in R_w$. But this implies that $w \in T_{p(v)}$ contradicting its definition. 
If $p(w) = p(v)-1$, then also $v$ and $w$ create a detour of remainder $2$ and, by the assumption on $v$, the two other vertices in $T$ creating sidewalks of remainder $1 \bmod 5$ are distinct from $w$. Thus, the detour can be combined with the two sidewalks and we obtain a contradiction. 

\subsection{Finalization in the cases with $\sum_i \mu_i + \sum_t \eps_t = 1$}

Note that $\sum_{i=1}^k \mu_i + \sum_{t\in T} \eps_t = 1$ implies that $\sum_{t\in T} |R_t| = k-2$ and $\ell \bmod k = 3$, and the family of non-empty sets $R_t$ is one of the six families listed in Table~\ref{tab:special_cases}. Intuitively, the sets of remainders are saturated, and any additional remainder leads to a contradiction. 

If $\sum_{i=1}^k \mu_i =1$, then not every vertex $t \in T$ having $1 \in R_t$ needs to have an outneighbor not following the cycle structure, so the conclusion of Observation~\ref{obs:special_cases_escape} does not hold. The following two observations will allow us to control the irregular blobs and vertices, and play the role of Observation~\ref{obs:special_cases_escape}.

\begin{obs}\label{obs:special_cases_1with_escape}
If $\sum_i \mu_i + \sum_t \eps_t = 1$ and a vertex $v$ satisfies $1 \in R_v$ and $\eps_v=0$, then every outneighbor $u$ of $v$ outside $V_{p(v)+1} \cup T_{p(v)}$ satisfies $u \in T$, $\eps_u = 1$ and $O_u = \{p(v)+2\}$.
\end{obs}

\begin{proof}
By Observation~\ref{obs:eps0}, $O_v = \{p(v)+1\}$, so all outneighbors of $v$ in $V' \setminus V_{p(v)+1}$ lie in $X$. Hence, $u \in X \cap V_j$ for some $j \neq p(v)+1$, or $u \in T \setminus T_{p(v)}$.

In the former case, $u \in X \cap V_j$ for some $j \neq p(v)+1$, and for each $i \in I_v$ there is a detour from $V_i$ through $v$ to $u$ of remainder $i-j+2$. If such a remainder is not in $R_v \cup \{0\}$, then we get a contradiction by Observation~\ref{obs:join_sidewalks}. Therefore, all the $|I_v|$ remainders $i-j+2$ are in $R_v \cup \{0\}$. Since $\eps_v=0$, we have $|R_v|+1=|I_v|$, so $I_v - j + 2 = R_v \cup \{0\}$. But also $I_v - (p(v)+1) + 2 = R_v \cup \{0\}$, so $R_v \cup \{0\}$ is invariant under the translation by $p(v)+1-j \neq 0$, which is impossible, as in the proof of Observation~\ref{obs:eps0}. Therefore, $u \in T \setminus T_{p(v)}$. In particular, $p(v) \not\in I_u$.

Let $j \in O_u$. Since $\eps_v=0$, for every $r \in R_v \cup \{0\}$ there is a detour from $V_{r+p(v)-1}$ through $vu$ to $V_{j}$ of remainder $r+p(v)-j+2$. We claim that $p(v)-j+2 \in R_u \cup \{0\}$. Otherwise, $p(v)-j+2 \neq 0$, and the detours for $r=0$ and $r=1$ have the consecutive remainders $p(v)-j+2$ and $p(v)-j+3$, none of which is $3$ by Observation~\ref{obs:join_sidewalks}. Hence, $p(v)-j+2=1$, or $k=5$ and $p(v)-j+2=4$. In the former case, $1 \notin R_u$ and the detour for $r=1$ has remainder $2$. If $2 \in R_v$, then the detour for $r=2$ has remainder~$3$. Otherwise, there is another vertex $v'$ with $1 \in R_{v'}$, distinct from $u$ as $1 \notin R_u$, and we combine the detour with a sidewalk of remainder $1$ centred at $v'$. In the latter case, $4 \notin R_u$ and the detour for $r=0$ has remainder $4$. If $4 \in R_v$, then the detour for $r=4$ has remainder~$3$. Otherwise, the vertex with remainder $4$ is distinct from $v$ and $u$, and we combine the detour with a sidewalk of remainder $4$ centred at it. In each case Observation~\ref{obs:join_sidewalks} gives a contradiction, which proves the claim that $p(v)-j+2 \in R_u \cup \{0\}$.

Since also $i-j+2 \in R_u \cup \{0\}$ for every $i \in I_u$, and $p(v) \notin I_u$, the $|I_u|+1$ distinct values $i-j+2$ with $i \in I_u \cup \{p(v)\}$ lie in $R_u \cup \{0\}$. Thus, $|I_u| \leq |R_u|$, so $\eps_u = 1$ and $(I_u \cup \{p(v)\}) - j + 2 = R_u \cup \{0\}$. As this holds for every $j \in O_u$ and $R_u \cup \{0\}$ is not invariant under a non-zero translation, $O_u = \{j\}$. Then $R_u = (I_u - j + 2) \setminus \{0\}$ has $|I_u|$ elements, so $I_u - j + 2 = R_u$. Hence, $p(v)-j+2 = 0$, that is, $O_u = \{p(v)+2\}$.
\end{proof}

\begin{obs}\label{obs:special_cases_escape1}
If $\sum_i \mu_i + \sum_t \eps_t = 1$ and a vertex $v$ satisfies $1 \in R_v$ and $\eps_v=0$, then $\mu_{p(v)}=1$ or $v$ has an outneighbor $u$ in $T \setminus T_{p(v)}$ that satisfies $\eps_u=1$ and $O_u = \{p(v)+2\}$.
\end{obs}

\begin{proof}
The assumptions $\eps_v=0$ and $1 \in R_v$ imply that $v \in T_{p(v)}$, so by 
\[d^+(v) \geq \delta^+(G) = |V_{p(v)+1}| +  |T_{p(v)}| - \mu_{p(v)} = |V_{p(v)+1}| +  |T_{p(v)}\setminus\{v\}| +1 - \mu_{p(v)},\]
either $\mu_{p(v)}=1$ or $v$ has an outneighbor $u$ outside $V_{p(v)+1} \cup T_{p(v)}$. In the latter, by Observation~\ref{obs:special_cases_1with_escape}, $u \in T \setminus T_{p(v)}$, $\eps_u = 1$ and $O_u = \{p(v)+2\}$.
\end{proof}

We are ready to finalize the proof in the remaining cases. 

\smallskip
\noindent\textbf{Cases with $\sum_i \mu_i = 1$:}
\smallskip

Since $\sum_i \mu_i = 1$, we have $\sum_{t \in T} \eps_t = 0$. 
All the considered families contain a vertex $v$ with $\{1,2\} \subseteq R_v$ or two vertices $u,v$ with $1 \in R_u \cap R_v$. 
In the former, Observation~\ref{obs:special_cases_escape1} gives $\mu_{p(v)} = 1$, while by Observation~\ref{obs:2sidewalk}, $\mu_{p(v)+1}=1$, which gives a contradiction with~\eqref{eq:totaldefect}. 
In the latter, by Observation~\ref{obs:special_cases_escape1}, $\mu_{p(v)} = \mu_{p(u)} = 1$, so $p(u)=p(v)$. Among vertices $v$ and $u$, there is at least one that does not have an arc to the other vertex. We may assume it is $v$. Since $u,v \in T_{p(v)}$, 
\[d^+(v) \geq \delta^+(G) = |V_{p(v)+1}| +  |T_{p(v)}| - \mu_{p(v)} = |V_{p(v)+1}| +  |T_{p(v)}\setminus\{u,v\}| + 1,\]
so $v$ has an outneighbor $w$ outside $V_{p(v)+1} \cup T_{p(v)}$. By Observation~\ref{obs:special_cases_1with_escape}, $\eps_w=1$, which is a contradiction. 

\smallskip
\noindent\textbf{Cases with $\sum_t \eps_t = 1$:}
\smallskip

Let $z \in T$ be the vertex with $\eps_z=1$. 
The assumption $\sum_t \eps_t = 1$ gives $\sum_i \mu_i = 0$, so no set $R_t$ contains $2$ by Observation~\ref{obs:2sidewalk}. Thus, there are two vertices $u,v\in T$ with $1 \in R_u \cap R_v$. 

Assume first that both $u$ and $v$ are distinct from $z$. Then $\eps_v = \eps_u = 0$ and by Observation~\ref{obs:special_cases_escape1}, $vz, uz \in E(G)$ and $O_z = \{p(v)+2\}=\{p(u)+2\}$, so $p(v)=p(u)$. As before, we may assume that $v$ has no arc to $u$, and since
\[d^+(v) \geq \delta^+(G) = |V_{p(v)+1}| +  |T_{p(v)}| - \mu_{p(v)} = |V_{p(v)+1}| +  |T_{p(v)}\setminus\{u,v\}| + 2,\]
vertex $v$ has at least two outneighbors outside $V_{p(v)+1} \cup T_{p(v)}$. The one distinct from $z$ gives a contradiction by Observation~\ref{obs:special_cases_1with_escape}.

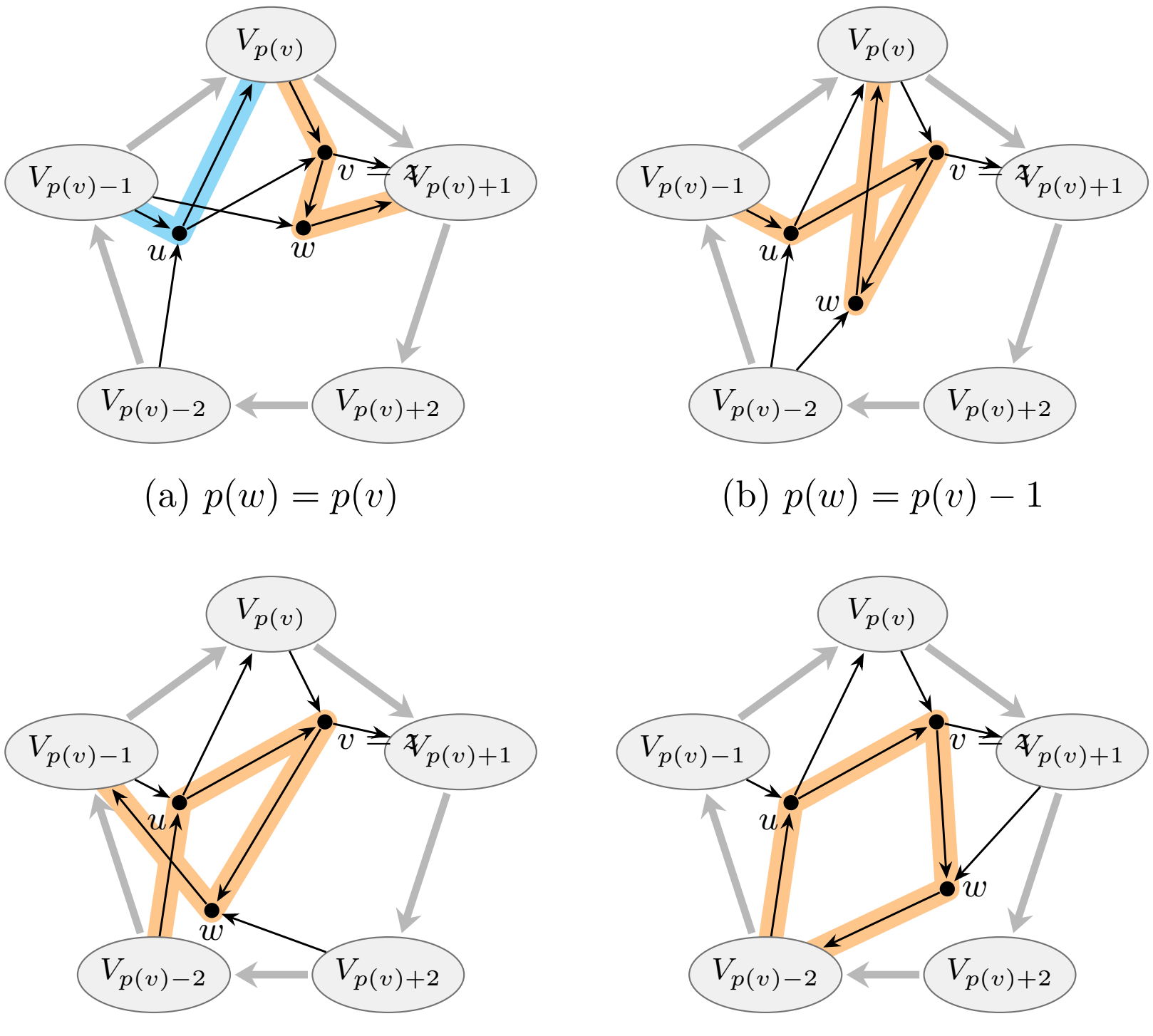
\begin{figure}[ht]
\centering
\begin{tikzpicture}
  \blobsfive
  \node[vtx,label={[vlab]below left:$u$}] (u) at (-0.85,0.1) {};
  \node[vtx,label={[vlab]below right:$v=z$}] (v) at (0.5,0.85) {};
  \node[vtx,label={[vlab]below:$w$}] (w) at (0.3,0.15) {};
  \begin{scope}[on background layer]
    \draw[hl] (B0.center) -- (v.center) -- (w.center) -- (B1.center);
    \draw[hlb] (Bm1.center) -- (u.center) -- (B0.center);
  \end{scope}
  \draw[arc] (Bm2) -- (u);
  \draw[arc] (Bm1) -- (u);
  \draw[arc] (u) -- (B0);
  \draw[arc] (B0) -- (v);
  \draw[arc] (v) -- (B1);
  \draw[arc] (Bm1) -- (w);
  \draw[arc] (w) -- (B1);
  \draw[arc] (u) -- (v);
  \draw[arc] (v) -- (w);
  \node[pcap] at (0,-2.35) {(a) $p(w)=p(v)$};
\end{tikzpicture}
\hspace{6mm}
\begin{tikzpicture}
  \blobsfive
  \node[vtx,label={[vlab]below left:$u$}] (u) at (-0.85,0.1) {};
  \node[vtx,label={[vlab]below right:$v=z$}] (v) at (0.5,0.85) {};
  \node[vtx,label={[vlab]left:$w$}] (w) at (-0.25,-0.55) {};
  \begin{scope}[on background layer]
    \draw[hl] (Bm1.center) -- (u.center) -- (v.center) -- (w.center) -- (B0.center);
  \end{scope}
  \draw[arc] (Bm2) -- (u);
  \draw[arc] (Bm1) -- (u);
  \draw[arc] (u) -- (B0);
  \draw[arc] (B0) -- (v);
  \draw[arc] (v) -- (B1);
  \draw[arc] (Bm2) -- (w);
  \draw[arc] (w) -- (B0);
  \draw[arc] (u) -- (v);
  \draw[arc] (v) -- (w);
  \node[pcap] at (0,-2.35) {(b) $p(w)=p(v)-1$};
\end{tikzpicture}
\\[4mm]
\begin{tikzpicture}
  \blobsfive
  \node[vtx,label={[vlab]below left:$u$}] (u) at (-0.85,0.1) {};
  \node[vtx,label={[vlab]below right:$v=z$}] (v) at (0.5,0.85) {};
  \node[vtx,label={[vlab]below:$w$}] (w) at (-0.55,-0.9) {};
  \begin{scope}[on background layer]
    \draw[hl] (Bm2.center) -- (u.center) -- (v.center) -- (w.center) -- (Bm1.center);
  \end{scope}
  \draw[arc] (Bm2) -- (u);
  \draw[arc] (Bm1) -- (u);
  \draw[arc] (u) -- (B0);
  \draw[arc] (B0) -- (v);
  \draw[arc] (v) -- (B1);
  \draw[arc] (B2) -- (w);
  \draw[arc] (w) -- (Bm1);
  \draw[arc] (u) -- (v);
  \draw[arc] (v) -- (w);
  \node[pcap] at (0,-2.35) {(c) $p(w)=p(v)-2$};
\end{tikzpicture}
\hspace{6mm}
\begin{tikzpicture}
  \blobsfive
  \node[vtx,label={[vlab]below left:$u$}] (u) at (-0.85,0.1) {};
  \node[vtx,label={[vlab]below right:$v=z$}] (v) at (0.5,0.85) {};
  \node[vtx,label={[vlab]right:$w$}] (w) at (0.6,-0.7) {};
  \begin{scope}[on background layer]
    \draw[hl] (Bm2.center) -- (u.center) -- (v.center) -- (w.center) -- (Bm2.center);
  \end{scope}
  \draw[arc] (Bm2) -- (u);
  \draw[arc] (Bm1) -- (u);
  \draw[arc] (u) -- (B0);
  \draw[arc] (B0) -- (v);
  \draw[arc] (v) -- (B1);
  \draw[arc] (B1) -- (w);
  \draw[arc] (w) -- (Bm2);
  \draw[arc] (u) -- (v);
  \draw[arc] (v) -- (w);
  \node[pcap] at (0,-2.35) {(d) $p(w)=p(v)+2$};
\end{tikzpicture}
\caption{Cases with $\sum_t \eps_t = 1$, drawn for $k=5$, with all possible positions of $w$.}\label{fig:eps1}
\end{figure}

Therefore, $z$ is one of the vertices $u$ and $v$, by symmetry, we may assume that $v=z$. 
By Observations~\ref{obs:special_cases_escape1} and \ref{obs:special_cases_1with_escape}, $uv \in E(G)$ and $O_v = \{p(u)+2\}$, so $p(v)=p(u)+1$. Since $v \in T_{p(v)}$ as $1 \in R_v$, by 
\[d^+(v) \geq \delta^+(G) = |V_{p(v)+1}| +  |T_{p(v)}| - \mu_{p(v)} = |V_{p(v)+1}| +  |T_{p(v)}\setminus\{v\}| + 1,\]
$v$ has an outneighbor $w$ outside $V_{p(v)+1} \cup T_{p(v)}$, so with $p(w) \neq p(v)+1$. 
The possible positions of $w$ are shown in Figure~\ref{fig:eps1}.
If $p(w)=p(v)$, then $v$ and $w$ create a detour of remainder $2$ that can be combined with a sidewalk of remainder $1$ centered at $u$, leading to a contradiction by Observation~\ref{obs:join_sidewalks}. 
While if $p(w)=p(v)-1$ or $p(w)=p(v)-2$, then $u$, $v$ and $w$ create a forbidden detour of remainder $3$.
This already gives a contradiction in the case $k=4$, while for $k=5$ we have $p(w)=p(v)+2$ and $u$, $v$, $w$ create a detour of remainder~$4$.
To avoid contradiction by Observation~\ref{obs:join_sidewalks}, in case of the family $(\{1\}, \{1\}, \{4\})$, we have $4 \in R_w$, which implies that $w \in T_{p(v)}$, a contradiction. 
Whereas, in case of the family $(\{1\}, \{1,4\})$, we must have $4 \in R_u$ or $4 \in R_v$. Then $u$, $v$, $w$ create a forbidden detour of remainder $3$, a contradiction.

\section*{AI statement}

All results in this paper have been obtained without AI assistance. The paper itself, except Figures 2--6, has also been written fully by the authors. The authors have used Claude Opus 5.5 only for linguistic checks, a final correctness check and making figures. 

\bibliography{refs}
\bibliographystyle{plain}

\end{document}